\documentclass[12pt]{extarticle}
\usepackage{amsmath, amsthm, amssymb, color}
\usepackage[colorlinks=true,,citecolor=red, linkcolor=blue,urlcolor=blue]{hyperref}
\usepackage{graphicx}
\usepackage{caption}
\usepackage{mathtools}
\usepackage{enumerate}
\usepackage{verbatim}
\usepackage{tikz,tikz-cd,tikz-3dplot}
\usepackage{amssymb}
\usetikzlibrary{matrix}
\usetikzlibrary{arrows}
\usepackage{algorithm}
\usepackage[noend]{algpseudocode}
\usepackage{caption}
\usepackage[normalem]{ulem}
\usepackage{subcaption}
\usepackage{multicol}
\usepackage{makecell}
\usepackage{array}
\usepackage{enumitem}

\oddsidemargin \evensidemargin
\newtheorem{theorem}{Theorem}[section]

\newtheorem{proposition}[theorem]{Proposition}
\newtheorem{lemma}[theorem]{Lemma}
\newtheorem{corollary}[theorem]{Corollary}

\newtheorem{conjecture}[theorem]{Conjecture}
\theoremstyle{definition}

\newtheorem{remark}[theorem]{Remark}
\newtheorem{cor}[theorem]{Corollary}
\newtheorem{example}[theorem]{Example}

\newcommand{\PP}{\mathbb{P}}
\newcommand{\RR}{\mathbb{R}}
\newcommand{\QQ}{\mathbb{Q}}
\newcommand{\CC}{\mathbb{C}}
\newcommand{\ZZ}{\mathbb{Z}}

\title{\bf Gram Matrices for Isotropic Vectors}
\author{Yassine El Maazouz, Bernd Sturmfels and Svala Sverrisd\'ottir}

\date{}

\begin{document}

\maketitle

\begin{abstract}
\noindent
We investigate determinantal varieties
for symmetric matrices that have
zero blocks along the main diagonal.
In theoretical physics,  these arise as Gram matrices for
kinematic variables in quantum field theories.
We study the ideals of relations among
functions in the matrix entries 
that serve as building blocks for conformal correlators.
 \end{abstract}

\section{Introduction}

Consider six row vectors $P_1,W_1,P_2,W_2,P_3,W_3$
in $\CC^r$ which are isotropic in the sense that
\begin{equation}
\label{eq:PWeqn}
 P_i \cdot P_i \,= \,P_i \cdot W_i \,= \,W_i \cdot W_i \,=\, 0 \quad
\hbox{for all $i$}. 
\end{equation}
The dot refers to some non-degenerate quadratic form on $\CC^r$.
The Gram matrix of pairwise inner products among the six vectors
has zero blocks of size  $2 \times 2$
along its main diagonal:
\begin{equation}
\label{eq:matrixX} \quad X \,\, = \,\, \begin{bmatrix}
  0 &  0 &  x_{13} & x_{14} & x_{15} & x_{16}  \\
  0 &  0  & x_{23} &  x_{24} &  x_{25} & x_{26} \\
  x_{13} & x_{23} & 0  & 0 &  x_{35} & x_{36} \\
  x_{14} & x_{24} & 0 &  0 &  x_{45} & x_{46} \\
  x_{15} & x_{25} & x_{35} & x_{45} & 0 &  0  \\
  x_{16} & x_{26} & x_{36} & x_{46} & 0 &  0 
  \end{bmatrix} 
 \,\,\,\, = \,\,\,\,
  \begin{bmatrix}   P_1 \\ W_1 \\ P_2 \\ W_2 \\ P_3 \\ W_3 \end{bmatrix} \cdot
   \begin{bmatrix}   P_1 \\ W_1 \\ P_2 \\ W_2 \\ P_3 \\ W_3 \end{bmatrix}^T .
\end{equation}
What are the constraints imposed on this Gram matrix
by the dimension $r$ of the space $\CC^r$?

\smallskip

Set-theoretically, the variety of such matrices is cut out by the
ideal of $(r+1) \times (r+1)$-minors of~$X$.
We now discuss this determinantal variety 
for each of the six possible values of~$r$.
The following list of six nested varieties will be examined  further in Example~\ref{ex:dreizwei}:
   \begin{itemize}
  \item[$r=0$:] The ideal contains all $12$ variables. Its variety is just the origin. The matrix is zero.
  \vspace{-0.2cm}
  \item[$r=1$:] The variety is the origin but with multiplicity $13$, given by
   the ideal of $2 \times 2$ minors.
  \vspace{-0.2cm}
\item[$r=2$:] The variety has
three irreducible components.
 There are also four embedded~primes.
   \vspace{-0.2cm}
 \item[$r=3$:] The variety is toric and hence irreducible.
The ideal is primary, but not radical.
  \vspace{-0.2cm}
\item[$r=4$:]   The ideal is radical but not prime.
It has one main component and three others.
  \vspace{-0.2cm}
\item[$r=5$:] The determinant of $X$ is irreducible.
It generates a principal ideal that is prime.
  \end{itemize}
  
Our study of such determinantal varieties is motived by theoretical physics, namely by
conformal field theory \cite{CPPR} and its recent use in cosmology \cite{BMPR}.
Here, the vectors $P_i$ and $W_i$ are momentum vectors and polarization vectors, respectively.
Physical quantities are functions of the Gram matrix $X = [x_{ij}]$.
These functions must satisfy the additional property of being invariant when
$W_i$ is replaced by $W_i + \alpha_i P_i$.
The following functions have this property:
\begin{equation}
\label{eq:fundinv}
\begin{matrix}
P_{ij} & = &  x_{2i-1,2j-1},\\
H_{ij} & = & x_{2i-1,2j-1} \cdot x_{2i,2j}\,-\, x_{2i,2j-1} \cdot x_{2j,2i-1},  \\
V_{i,jk} & = &  \big(x_{2i,2j-1} \cdot x_{2k-1,2i-1} \,-\, x_{2i,2k-1}\cdot x_{2j-1,2i-1}\big)/x_{2j-1,2k-1}.
\end{matrix}
\end{equation}
These are the
{\em conformally-invariant structures} which were derived
in \cite[Section 4.2.2]{CPPR}:
\begin{equation}
\label{eq:structures}
 \begin{matrix}
P_{ij}  &=& P_i \cdot P_j, \\
H_{ij} &=& (W_i \cdot W_j) (P_i \cdot P_j) - (W_i \cdot P_j)(W_j \cdot P_i), \\
 V_{i,jk} &= &\bigl[(W_i \cdot P_j) (P_k \cdot P_i) - (W_i \cdot P_k)(P_j \cdot P_i)\bigr]/({P_j \cdot P_k}).
\end{matrix}
\end{equation}

\begin{remark}
    The numerator in the expression of $V_{i, jk}$ is already invariant under the translations $W_i \mapsto W_i + \alpha_i P_i$. However, the quantities of interest in conformal field theory are the \emph{three-point functions} which are polynomials in the $H_{ij}$ and $V_{i,jk}$. See \cite[Equation 4.17]{CPPR}.
\end{remark}
The expressions in \eqref{eq:fundinv} map the determinantal varieties of $X$ to different varieties in projective space.
We seek to find the ideal in the ring $\CC[P_{ij},H_{ij},V_{i,jk}]$ which describes the image of this map.
For instance, when $r = 5$, this ideal is generated by the expression
\begin{equation}
\label{eq:3gen} {\rm det}(X)\,\, = \,\,
4H_{12}H_{23}H_{13} - (V_{1,23}H_{23} - V_{2,13}H_{13} + V_{3,12}H_{12} + V_{1,23}V_{2,13}V_{3,12})^2,
\end{equation}
which appears in \cite[equation (2.10)]{BMPR} and \cite[equation (4.48)]{CPPR}.
Additionally, when $r=4$, the main component of the image variety is cut out by the following ideal
\begin{equation}
\label{eq:threebinomials} \bigl\langle \,
H_{12} + V_{123} V_{213}\,,\,
H_{13} - V_{123} V_{312}\,,\,
H_{23} + V_{213} V_{312} \ \bigr\rangle. 
\end{equation}
Our guiding question is how 
 to encode low rank Gram matrices 
like (\ref{eq:matrixX}) by ideals like~(\ref{eq:threebinomials}).

We now discuss the organization and main results of this paper.
In Section \ref{sec2} we study the variety of
symmetric $kn \times kn$ matrices
of rank $r$ with zero diagonal blocks of size $k \times k$.
Theorem \ref{thm:Vmain} states that the ideal of
$(r+1) \times (r+1)$ minors is prime for $r \geq 2k+1$.
The proof is based on a parametrization of our variety which rests on the
orthogonal Grassmannian of $k$-planes in $\CC^r$.
This variety has gained prominence in
cosmology  through the work in \cite{ABLPR}.

In Section \ref{sec3} we turn to 
theoretical physics. The embedding space formalism in
conformal field theory takes momentum vectors in $\RR^d$
onto the Lorentz cone in $\RR^r$ where $r=d\!+\!2$.
Conformal symmetries in $\RR^d$ are thus encoded by
the action of  the Lorentz group ${\rm SO}(1,r-1)$ on $\RR^r$.
We consider $n$ fields, given by
momentum vectors $P_i$ and polarization vectors $W_i$
for $i=1,2,\ldots,n$. This leads to $2n \times 2n$ Gram matrices of rank $r$ with
blocks of size~$k=2$.

Section \ref{sec4} addresses the case $d=2$.
Spinor-helicity formalism is used in Proposition~\ref{prop:spinorhelicity}
to parametrize the main component of $\mathcal{V}_{2,n,4}$.
We study the ideal of relations
among the conformally invariant structures
(\ref{eq:structures}).
In Theorem~\ref{thm:twopointvariety}, we present a
dimension formula for the varieties of  two-point functions.
This is proved using methods from tropical geometry~\cite{MS}.
Section \ref{sec5} concerns the same ideals  for $d=3$.
This is the case of interest for inflationary cosmology~\cite{BMPR}.
Implicitization is difficult, but we show the
state of the art for $n=4$ fields.
Our results are the foundation for current advances~\cite{BKMP}
on constraining conformal correlators.

Section \ref{sec6} departs from physics and returns us to
commutative algebra.
 We investigate the variety of
 Gram matrices with zero blocks in
  the low rank regime $r \leq 2k$. 
  Our main results (Theorems \ref{thm:dimension}
  and \ref{thm:components}) characterize
  all irreducible components 
  and their dimensions.
  
In conclusion, this article significantly strengthens the connection between theoretical physics and commutative algebra.
For physicists, it offers new tools for designing bootstrap methods for scattering amplitudes and conformal correlators.
Algebraists will find sharp theorems on a natural class of
determinantal varieties, leading to interesting open problems.

\section{Determinantal Varieties}
\label{sec2}

We fix a symmetric $kn \times kn$ matrix $X = (x_{ij})$
with indeterminate entries subject to the constraint that the
$k \times k$ blocks along the diagonal are zero. In symbols, the entries satisfy
\begin{equation}
\label{eq:postulate}
 \begin{matrix} & \qquad \qquad x_{ij} \,=\, x_{ji} & \!\! \hbox{for $\,\,1 \leq i \leq j \leq kn$ } & \\
{\rm and} & \,\, x_{jk+l,jk+m}\,  =\,	0 &\,\,
\hbox{for $\,j=0,1,\ldots,n-1$} & \hbox{and $\,\,1 \leq l \leq m \leq k$. }\end{matrix} 
\end{equation}
Thus the number of distinct entries $x_{ij}$ in the $kn \times kn$ matrix $X$ is equal to
\begin{equation}
\label{eq:Mvariables}  M \,\,= \,\,\binom{kn+1}{2} \,-\, n\binom{k+1}{2}. 
\end{equation}
We write $\CC[X]$ for the polynomial ring in these $M$ variables, and
we consider the affine variety 
$\mathcal{V}_{k,n,r}$ in the ambient affine space $\CC^M$
 whose points are the matrices $X$ of rank $\leq r$.

 In the special case $n=2$, there is only one block, namely
$X = \begin{scriptstyle} \begin{bmatrix} 0 \! & \! A \\ A^T \!\!\! &\! 0 \end{bmatrix} \end{scriptstyle} $
for some $A \in \CC^{k \times k}$.

\begin{proposition}
The variety $\,\mathcal{V}_{k,2,r}$ is
irreducible. To be precise, we have
\begin{equation} \label{eq:k2r}
\mathcal{V}_{k,2,r} \,=\, \{ \,A \in \CC^{k \times k} :
{\rm rank}(A) \leq \ell \,\} \quad \hbox{where
$\,\ell =  \lfloor r/2 \rfloor $.}
\end{equation}
The prime ideal given by the
 $(\ell +1) \times (\ell+ 1)$ minors of $A$
is the radical of the ideal of $(r+1) \times (r+1)$ minors of $X$.
The latter
has the ideal of $\ell \times \ell$ minors of $A$ as an embedded~prime.
\end{proposition}
 
 \begin{proof}
 Since ${\rm rank}(X) = 2 \cdot {\rm rank}(A)$, we have
  ${\rm rank}(X) \leq r$ if and only if
 ${\rm rank}(A) \leq \ell$. The last claim holds since
   $I_{r+1}(X):I_{\ell+i}(A) = I_{\ell}(A)$
for $i=2$ if $r$ is odd
and $i=1$ if $r$ is even.
 \end{proof}

Another special case for immediate deliberation
is that of blocks having size $1 \times 1$.

\begin{example}[$k=1$] \label{ex:k1}
Here $X$ is a symmetric $n \times n$ matrix with zeros on the diagonal.
The variety $\mathcal{V}_{1,n,r}$ is irreducible for $r \geq 3$. This was shown
for $r=3$ in \cite[Theorem 3.5]{DFRS}, where $\mathcal{V}_{1,n,3}$ was
identified with the {\em squared Grassmannian} ${\rm sGr}(2,n)$. The proof
for $r \geq 4$ is analogous to that given in \cite[Section 3]{DFRS}. 
For $r =0$, our variety is just the origin.
For $r=1$, the variety is still the origin, but  now
the $2 \times 2$-minors define the origin with multiplicity $2^{n-1}$.
The most interesting case arises for $r=2$.  This is the content of the next proposition.
\end{example}

\begin{proposition}\label{prop:k=1decomp}
The $3 \times 3$ minors of $X$ in Example \ref{ex:k1}
generate an unmixed radical ideal of codimension $\binom{n-1}{2}$
and degree $\binom{2n-3}{n-2}$. The variety $\,\mathcal{V}_{1,n,2}$
has $2^{n-1} -1 $ irreducible components.
\end{proposition}

\begin{proof}
For any partition  of $[n] = \{1,2,\ldots,n\}$
into two non-empty subsets $I$ and $J$,
we shall identify an irreducible subvariety of
codimension $\binom{n-1}{2}$ and degree $\binom{n-2}{| I | - 1}$.
The number of such partitions is
$2^{n-1}-1$, and the degrees add up to $\binom{2n-3}{n-2}$.
The component consists of all matrices $X$ such that
the submatrix $X_{I,J}$ with row indices $I$ and column indices $J$
has rank $\leq 1$ and $x_{rs} = 0$ whenever $\{r,s\}$ is contained either in $I$ or in $J$. This indeed has~codimension 
$$ \begin{matrix} (| I | - 1)(|J|-1) + \binom{|I|}{2} + \binom{|J|}{2} \,\,=\,\, \binom{n-1}{2}. \end{matrix} $$

The codimension and degree in our assertion matches that for the
$3 \times 3$ minors of  a generic symmetric $n \times n$  matrix,
by \cite[Proposition~12]{HT}. The diagonal entries form
a regular sequence modulo that Cohen-Macaulay ideal.
This is proved as in \cite[Section 3]{DFRS}, and we will
derive a more general result in Section \ref{sec6}.
This shows that each of the irreducible varieties above
is a component, these are all the components, and there are no
embedded primes.
\end{proof}

\begin{remark} Conca \cite{Conca} proved that the $(r+1) \times (r+1)$
minors of a generic symmetric matrix form a Gr\"obner basis.
This is false in our setting. Specifically, for
$k=1,n=4,r=2$, the ideal is generated by the
ten $3 \times 3$ minors of $X$. These do not form
a Gr\"obner basis for any term order.
One can show that any Gr\"obner basis
 must contain elements of degree four.
 \end{remark}

We now come to our main result in this section. 
It concerns the case of large rank $r$.
 
\begin{theorem} \label{thm:Vmain}
Let $r \geq 2k+1$. The determinantal variety $\mathcal{V}_{k,n,r}$ is irreducible
of dimension
\begin{equation}
\label{eq:expecteddimension}
    \dim(\mathcal{V}_{k,n,r}) \,\,=\,\, n\left(kr - \binom{k+1}{2}\right) - \binom{r}{2} .
\end{equation}
    The $(r+1) \times (r+1)$ minors of $\,X$ generate the prime ideal of $\,\mathcal{V}_{k,n,r}$,
    which is Cohen-Macaulay.
\end{theorem}

\begin{proof}
    We consider the variety of all symmetric $kn \times kn$ matrices of
    rank $\leq r$. It lives in $\CC^{\binom{kn+1}{2}}$ and spans this space.
    Every matrix in this variety can be factored as 
    $YY^T$ where $Y$ is a  $kn \times r$ matrix.
    We think of $Y$ as $n$ generic $k \times r$ matrices $Y_1, \dots, Y_n$ stacked one on top of the other, and we write $\CC[Y]$ for the corresponding polynomial ring.
    The coordinate ring of our variety is the subalgebra
    generated by the entries of the $kn \times kn$ matrix $YY^T$:
    $$ R \,\, = \CC[YY^T] \,\,= \,\, \CC[\,Y_i Y_j^T \,:\, 1 \leq i \leq j \leq n \,]  \quad \subset \,\,\, \CC[Y]. $$
    Note that $R$ is a Cohen-Macaulay domain because
    it coincides with the invariant ring  $\CC[Y]^{{\rm O}(r)}$
    of the action of the orthogonal group ${\rm O}(r)$  on $\CC^{kn \times r}$; see e.g.~\cite{DCP}.
    We consider the quotient of $R$ by the $n \binom{k}{2}$ quadrics that correspond to the diagonal block entries (\ref{eq:postulate}) of $X = Y Y^T$.
    
    That quotient is the coordinate ring of
    the determinantal variety $\mathcal{V}_{k,n,r}$  we wish to analyze:
    \begin{equation}
    \label{eq:CCV} \CC[\mathcal{V}_{k,n,r}] \,\, = \,\, R\, /\langle\, Y_i Y_i^T \,: \, i =1,\ldots,n \,\rangle. 
    \end{equation} 
    
    Since ${\rm O}(r)$ is a reductive algebraic group, its invariant ring $R := \CC[Y]^{{\rm O}(r)}$ 
    is  also Cohen-Macaulay by the Hochster-Roberts Theorem. It is hence    is a direct summand of $\CC[Y]$. 
    More precisely,
    using the Reynolds operator   $\rho \colon \CC[Y] \longrightarrow R = \CC[Y]^{{\rm O}(r)}$, we 
    can write 
    \[  \CC[Y] \,\,= \,\,R \oplus {\rm ker}(\rho) \quad \text{as } \text{$R$-modules}.
    \] Therefore, to show that the ring $\CC[\mathcal{V}_{k,n,r}]$ is an integral domain, it suffices to show that the ideal generated by the entries of the $n$ matrices $ Y_i Y_i^T$ is prime in $\CC[Y]$.
    
    Suppose that $r \geq 2k+1$. We claim that the following 
    $\CC$-algebra is an integral~domain:
    $$  
        \CC[Y]  \, /\langle\, Y_i Y_i^T \,: \, i =1,\ldots,n \,\rangle \,= \,
        \CC[Y_1] \, /\langle Y_1 Y_1^T \rangle \otimes_\CC 
        \CC[Y_2] \, /\langle Y_2 Y_2^T \rangle \otimes_\CC  \cdots \otimes_\CC
        \CC[Y_n] \, /\langle Y_n Y_n^T \rangle.
    $$
    This tensor product decomposition is as in the proof of \cite[Theorem 2.1]{RSS}.
    The factors are all isomorphic. We now write $Y = (y_{ij})$ for a generic $k \times r$ matrix.
    Since the tensor product of integral domains over $\CC$ is a domain, we must
    show that the following ideal is prime:
    \begin{equation}
    \label{eq:YYT}
     I \,\, = \,\,\bigl\langle \,Y  Y^T \,\bigr\rangle \,\, = \,\,
    \biggl\langle \,
    \sum_{\ell=1}^{r}
     y_{i \ell}\, y_{j \ell} : 1 \leq i \leq j \leq k  \biggr\rangle
     \,\,\, \subset \,\,\, \CC[Y].
    \end{equation}

    Let $T := \CC[Y] / I$  and write $\mathcal{Y} := {\rm Spec}(T)$. We consider the vector subbundle  $\mathcal{E}$ of the trivial vector bundle $\CC^{k \times r} \times \mathrm{OGr}(k,r)$ defined as follows:
    \[
        \mathcal{E} := \Big\{ (Y,U) \in \CC^{k \times r} \times \mathrm{OGr}(k,r) \colon {\rm rowspan}(Y) \subset U \Big\}.
    \]
    Here, ${\rm OGr}(k,r)$ is the orthogonal Grassmannian of $k$-dimensional isotropic subspaces of $\CC^r$. We write $\pi_1 \colon \mathcal{E} \to \CC^{k \times r}$ and $\pi_2  \colon \mathcal{E} \to {\rm OGr}(k,r)$ for the projections onto the first and second factors respectively. Since $r > 2k+1$, the Grassmannian ${\rm OGr}(k,r)$ is irreducible and $\mathcal{Y}(\CC)$ is exactly the image of $\pi_1$. So the variety $\mathcal{Y}(\CC)$ is also irreducible. Moreover $\pi_1 \colon \mathcal{E} \to \mathcal{Y}(\CC)$ is birational and surjective, so $\dim(\mathcal{Y}) = \dim(\mathcal{E})$. The projection map $\pi_2$ is also surjective and all of its fibers have dimension $k^2$. Hence $\dim(\mathcal{Y}) = \dim(\mathcal{E}) = k^2 + \dim({\rm OGr}(k,r))$.
    It is known that $\dim({\rm OGr}(k,r)) = k(r-k) - \binom{k+1}{2}$, see \cite[Proposition 2.4]{EM}. So we deduce that 
    \[
        \dim(\mathcal{Y}) = kr - \binom{k+1}{2}.
    \]
    Since the ideal $I$ is generated by $\binom{k+1}{2}$ elements, then $T = \CC[Y] / I$ is a complete intersection ring, and hence it is Cohen-Macaulay. In particular, $T$ satisfies Serre's depth condition ($\mathrm{S}_2$).
    
    We next claim that the ring $T := \CC[Y]/I$ is regular in codimension $1$.
    Consider the closed subscheme $\mathcal{Z}$ of $\mathcal{Y}$ which represents matrices of rank less than $k$.
    In symbols,
    \[
        \mathcal{Z} \,\,:=\,\, {\rm Spec} \Big( \CC[Y] \Big/ \big( I +  \langle  k \times k \text{ minors of } Y\rangle\big)    \Big).
    \]
    This closed subscheme has dimension $(k-1)r - \binom{k}{2}  + (k-1)$. Hence $\mathcal{Z}$ has codimension $r-2k+1 \geq 2$. Therefore, to prove our claim, it suffices to show that every closed point in the open set $\mathcal{U} := \mathcal{Y} \backslash  \mathcal{Z}$ is smooth. For this we first note that ${\rm GL}(k) \times \mathrm{O}(r)$ acts on $\CC[Y]$ via
    \[ 
        (A,B) \cdot Y \, =\, A\, Y B \quad \text{for } A  \in {\rm GL}(k) \,\text{ and } \,B \in \mathrm{O}(r).
    \]
    The ideal $I$ is invariant under this action. So ${\rm GL}(k) \times \mathrm{O}(r)$ acts on the scheme $\mathcal{Y}$. This action is transitive on $\mathcal{U}(\CC)$ because ${\rm GL}(k)$ acts transitively on $k \times r$ matrices of rank $k$ with the same row span,
    and the orthogonal group $\mathrm{O}(r)$ acts transitively on the orthogonal Grassmannian $\mathrm{OGr}(k,r)$. Indeed, this Grassmannian can be constructed as $\mathrm{O}(r)$ modulo a parabolic subgroup. 
    
    We now show that one, and hence every, closed point in $\mathcal{U}$ is smooth. We pick the point 
    \[
        Y_0 := \begin{bmatrix} I_k & i I_k & 0_{k \times (r-2k) }\end{bmatrix} \in \mathcal{U}(\CC).
    \]
    Let $A,B \in \CC^{k \times k}$ and $C \in \CC^{k \times (r-2k)}$ and set $X = \begin{bmatrix} A & B & C\end{bmatrix} \in \CC^{k \times r}$. We have
    \[
        (Y_0 + \epsilon X)(Y_0 + \epsilon X)^T = 0 \mod \epsilon^2 \quad \text{if and only if} \quad A + iB \text{ is skew-symmetric}.
    \]
    The map $\CC^{k \times k} \oplus \CC^{k \times k} \to \CC^{k \times k}, (A,B) \mapsto A + iB$ is surjective and its kernel has dimension $k^2$, and the space of skew-symmetric matrices has dimension $\binom{k}{2}$. 
    We conclude that the Zariski tangent space of $\mathcal{Y}$ at the point $Y_0$ has dimension
    \[
        \dim(T_{Y_0}\mathcal{U}) \,\,=\,\, k(r - 2k) + \binom{k}{2} + k^2 \,\,=\,\, kr - \binom{k+1}{2} \,\,=\,\, \dim(\mathcal{Y}). 
    \]
    
    So the point $Y_0 \in \mathcal{U}(\CC)$, and hence every point in $\mathcal{U}(\CC)$, is smooth. Since the singular locus of $\mathcal{U}$ is closed, and since any nonempty closed subset of the finite-type $\CC$-scheme $\mathcal{U}$ contains a closed point, it follows that $\mathcal{U}$ is a smooth scheme.
    This proves our claim which states that $T = \CC[Y]/I$ is regular in codimension $1$. In particular, the ring $T$ satisfies Serre's regularity condition $  (\mathrm{R}_1)$
    
    Using Serre's criterion for normality, we deduce that $T = \CC[Y]/I$ is a normal Noetherian ring. Hence $T = T_1 \times \dots \times T_m$ where each $T_i$ is a normal domain. Since $\mathcal{Y}$ is irreducible we must have $m=1$ and hence $T$ is a normal domain. We then conclude that $I$ is indeed a prime ideal.
    
    We have now shown that the $(r+1) \times (r+1)$ minors of $\,X$ generate
    the prime ideal of $\mathcal{V}_{k,n,r}$. The dimension formula in
    (\ref{eq:expecteddimension}) follows from our construction
    of $\mathcal{V}_{k,n,r}$ from $n$ copies of ${\rm Spec}(T)$,
    which has dimension $kr - \binom{k+1}{2}$,
    modulo a faithful action of ${\rm O}(r)$, which has dimension $\binom{r}{2}$.
To verify the Cohen-Macaulay property, we consider a
symmetric matrix of size $k n \times kn$ where all
$\binom{kn+1}{2}$ entries are distinct unknowns.
The ideal of $(r+1) \times (r+1)$ minors of such a matrix
 is prime and Cohen-Macaulay,
of codimension $\binom{kn-r+1}{2}$ and degree
$\prod_{j=0}^{kn-r-1} \binom{kn+j}{kn-r-j} / \binom{2j+1}{j}$.
See \cite[Proposition~12]{HT}.
If we add the expected  $\binom{kn-r+1}{2}$ to the
right hand side of  (\ref{eq:expecteddimension}) then
we obtain the ambient dimension $M$ in 
(\ref{eq:Mvariables}). This shows that the entries on
the $k \times k$ diagonal blocks form a regular sequence,
and hence the $(r+1) \times (r+1)$ minors of our matrix $X$
generate a Cohen-Macaulay ideal of the expected dimension.
\end{proof}

\begin{corollary} \label{cor:expecteddegree}
 If $r \geq 2k+1$ then the degree of our irreducible
determinantal variety equals
\begin{equation}
\label{eq:expecteddegree}
 {\rm degree}(\mathcal{V}_{k,n,r}) \,\,\,\, = \,\,
\prod_{j=0}^{kn-r-1} \binom{kn+j}{kn-r-j} / \binom{2j+1}{j}.
\end{equation}
\end{corollary}

\begin{remark} \label{rmk:severalsteps}
Several of the steps in the proof above remain valid for $r = 2k$. 
The Cohen-Macaulay property persists, and the
formulas   (\ref{eq:expecteddimension}) and
(\ref{eq:expecteddegree}) remain valid.
Inside the proof, the vector bundle $\mathcal{E}$ now has two distinct irreducible connected components. Consequently the scheme $\mathcal{Y} := {\rm Spec}(T)$ has two irreducible components and the same argument we used in Theorem 2.5 shows that each of these components has dimension $2k^2 - \binom{k+1}{2}$. So $T := \CC[Y]/I$ is a complete intersection ring and hence is still Cohen-Macaulay and satisfies Serre's depth condition (S$_2$). However, the locus $\mathcal{Z}$ now has codimension $1$ in ${\rm Spec}(T)$. Hence $T$ is singular in codimension~$1$, and 
 Serre's regularity condition $(\mathrm{R}_1)$ fails.
 In fact, the variety $\mathcal{Y}(\CC)$ is the union of two irreducible components.
 After a linear change of coordinates on the row space $\CC^r$, so that the inner product is the bilinear form $(u,v) \mapsto \sum_{i=1}^{r} (-1)^{i} u_i v_{i}$, these two components can be described by the following linear identities among $k \times k$ minors of the matrix $Y$:
\begin{equation}\label{eq:twocomponents} 
\hbox{
${\rm det}(Y_A) = 
{\rm det}(Y_{A^c})$
for all $A \in \binom{[2k]}{k}$ 
\quad or \quad
${\rm det}(Y_A) = -
{\rm det}(Y_{A^c})$
for all $A \in \binom{[2k]}{k}$.
}
\end{equation}
So, in the smooth locus of $\mathcal{Y}$ we see two isomorphic connected~components.
\end{remark}

\begin{example}[$k=2,r=3$]
The ideal $I = \langle Y Y^T \rangle$ is the complete intersection of
three quadrics in six unknowns. This ideal is primary but not prime.
In particular, every point of ${\rm Spec}(T)$
 is singular.
The prime ideal $\sqrt{I}$
is generated by $I$ and the three $2 \times 2$ minors of $Y$.
\end{example}

\begin{example}[$k=2,r=4$]
The ideal $I = \langle Y Y^T \rangle$ is generated by three quadrics in eight unknowns.
It is the intersection of two prime ideals of codimension three,
each generated by six quadrics. The three additional generators are the
determinant equations in (\ref{eq:twocomponents}).
\end{example}

We conclude this section by revisiting the rank strata we encountered in the Introduction.
The findings  in the following example will be explained by the results in Theorems
\ref{thm:dimension} and~\ref{thm:components}.

\begin{example}[$n=3,k=2$]  \label{ex:dreizwei}
 The $6 \times 6$ matrix $X$ in  (\ref{eq:matrixX}) contains $12$ distinct variables $x_{ij}$.
  The ideal defining  $\mathcal{V}_{2,3,r} = \{ X : {\rm rank}(X) \leq r \}$ in $\CC^{12}$ is generated by the
  $(r\!+\!1) \times (r\!+\!1)$ minors of $X$. We examine the six determinantal
  ideals that arise for the various ranks.
    \begin{itemize}
  \item[$r=0$:] The variety $\mathcal{V}_{2,3,0}$ is just the origin, so it has 
    codimension $12$. The point is reduced because the ideal is the maximal ideal
    which is
    generated by the $12$ unknowns $x_{ij}$.
\item[$r=1$:] The variety $\mathcal{V}_{2,3,1}$ is also just the origin, but now with multiplicity
$13$. The ideal is generated by the $2 \times 2$ minors of $X$, and it equals the square of the maximal ideal.
\item[$r=2$:] The variety $\mathcal{V}_{2,3,2}$ has codimension $7$ and degree $12$, with
three irreducible components, where one $2 \times 2$ block is zero
and the remaining $2 \times 4$ matrix has rank $1$.  The
ideal is generated by $145$ cubics, namely the $3 \times 3$ minors.
It is not radical. Besides the three minimal primes, there are four embedded primes,
given by intersecting the components.
\item[$r=3$:] The variety $\mathcal{V}_{2,3,3}$ is irreducible but not reduced. It is
a toric variety of codimension $6$ with a multiplicity $8$ structure. See (\ref{eq:hadamard}).
The ideal
of $4 \times 4$ minors of $X$ is generated~by $105$ quartics.
It is primary and has degree $112$.
Its radical has degree $14$. It is generated by the $15$ distinct 
$ 2 \times 2$ minors  of the $2 \times 4$ matrices formed by pairs of blocks. 
\item[$r=4$:]   The ideal of $5 \times 5$ minors has $21$ quintic generators. It is radical
of pure codimension $3$ and degree $35$. The variety $\mathcal{V}_{2,3,4}$ has four
irreducible components, of degrees $11 ,8,8,8$. The last three components are
extraneous in the sense that two of the three $2 \times 2$ blocks are singular on each component.
The prime ideal of the main component
$\widetilde{\mathcal{V}}_{2,3,4}$
 of degree $11$ is generated by nine cubics
such as $\,x_{13} x_{15} x_{46}-x_{13} x_{16} x_{45}-x_{14} x_{15} x_{36}+x_{14} x_{16} x_{35}$.
\item[$r=5$:] The determinant of $X$ is an irreducible polynomial of degree six
with $40$ terms. Its Jacobian ideal is radical, and it is the
intersection of seven prime ideals.
 The singular locus of the irreducible hypersurface
$\mathcal{V}_{2,3,5}$ strictly contains $\mathcal{V}_{2,3,4}$.
The latter contributes four of the seven components of the singular locus of $\mathcal{V}_{2,3,5}$.
The remaining three components are determinantal varieties,
given by the $2 \times 4$ blocks in $X$ having rank~$1$.
  \end{itemize}
The formula (\ref{eq:expecteddegree})
explains the degrees  
$112,35,6$ we computed
for $r=3,4,5$.
See also Table~\ref{tab:varieties1}.
  \end{example}

\section{From Cosmology to Nonlinear Algebra}
\label{sec3}

We now turn to connections with physics.
In conformal field theory \cite{CPPR}, one has $k=2$, and
there are $n$ fields or particles, which live in spacetime of dimension $d$.
The rank $r$ of our matrix $X$ satisfies $r=d+2$. This is the dimension of the
embedding space for the conformal group  in dimension $d$.
In the setting of cosmology  \cite{BMPR}, the parameters are $d = 3$ and  $r=5$.

Our article has its genesis in discussions with 
the team of Daniel Baumann, who works in theoretical cosmology \cite{ABLPR}.
In that setting, one considers $n$ fields on $\RR^3$, where the time is fixed, say 
 at the beginning of the hot big bang or the end of inflation.
 We use the   notation from~\cite{BMPR}.
The $i$th field is represented by a point $P_i = (p_{i0}:p_{i1}:p_{i2}:p_{i3}:p_{i4})$ in
the projective space $\PP^4$. This point lies on the {\em projective null cone}, shown in
 \cite[Figure 2]{BMPR}, so it satisfies
\begin{equation} \label{eq:quadricQ1} P_i \cdot P_i \,\, = \,\,
 p_{i0}^2 - p_{i1}^2 - p_{i2}^2 + p_{i3}^2 - p_{i4}^2 \,\, = \,\, 0 . 
 \end{equation}
We write $Q \subset \PP^4$ for the threefold
defined by this quadric. This is invariant under the
group ${\rm SO}(3,2)$, which plays the
role of the Lorentz group. We map $P_i$ to the skew-symmetric~matrix
\begin{equation} \label{eq:quadricQ2} {\bf P}_i \,\, = \,\,
\begin{bmatrix}
   0 & -p_{i3}-p_{i4} &     p_{i2} &  p_{i0}+p_{i1}  \\
   p_{i3}+p_{i4} &    0  &  p_{i1}-p_{i0} &   -p_{i2}  \\
  -p_{i2} & p_{i0}-p_{i1} &      0 &  p_{i3}-p_{i4}  \\
 -p_{i0}-p_{i1} &   p_{i2}  & p_{i4}-p_{i3} &    0   \\
\end{bmatrix}.
\end{equation}
The quadric in (\ref{eq:quadricQ1}) equals the {\em Pfaffian} of ${\bf P}_i$.
So, the determinant of ${\bf P}_i$ is the square of~(\ref{eq:quadricQ1}).

In the particle physics setting of \cite{PRRVZ, RSS},
the matrix (\ref{eq:quadricQ2}) plays the role of the
momentum space Dirac matrix.
Our aim is to introduce algebraic varieties that are
analogous to the varieties studied in \cite[Sections 5 and 6]{RSS},
but now in the cosmology context developed in~\cite{BMPR}.
The material that follows can be useful for 
bootstrap approaches to conformal correlators. These rest on
algebraic and combinatorial constraints that are encoded in our ideals;
see~\cite{BKMP}.

To this end we introduce $n$ further points
 $W_1,\ldots,W_n$ on the quadric $Q$.
 These are {\em null polarization vectors}.
 This means that the translation
 $P_i+W_i$ also lies on $Q$ for $i=1,\ldots,n$.
Thus the system (\ref{eq:PWeqn})
of $3n$ quadratic equations is assumed to hold
for our $2n$ vectors.
Now, the dot product refers to the bilinear form that is induced by the
quadratic form  (\ref{eq:quadricQ1}), i.e.
$$ P_i \cdot  W_i \,\, = \,\,
p_{i0} w_{i0}  -p_{i1} w_{i1}  -p_{i2} w_{i2}  +p_{i3} w_{i3}  -p_{i4} w_{i4} .$$
Choosing $n$  momenta $P_i$ and $n$ polarization vectors $W_i$ specifies 
a point in $(\PP^4)^n \times (\PP^4)^n = (\PP^4)^{2n}$. Equations (\ref{eq:PWeqn})
define a subvariety of codimension $3n$ in the ambient space $(\PP^4)^{2n}$.

\begin{proposition} \label{prop:abcdefg}
The subvariety of $ (\PP^4)^{2n}$ defined by  (\ref{eq:PWeqn}) has the  parametric representation
\begin{equation}\label{eq:Knpara} 
\resizebox{.92\hsize}{!}{$
\begin{aligned}
    p_{i0} &= a_i b_i d_i+a_i c_i e_i - b_i e_i + c_i d_i, \qquad &
    p_{i1} &= -a_i b_i d_i - a_i c_i e_i + b_i e_i + c_i d_i, \qquad &  p_{i2} &= -2 c_i e_i , \\ 
    p_{i3} &= b_i c_i d_i+c_i^2 e_i - a_i d_i + e_i, \qquad &
    p_{i4} &= -b_i c_i d_i-c_i^2 e_i-a_i d_i+e_i, & \\
    w_{i0} &= a_i b_i f_i + a_i c_i g_i - b_i g_i + c_i f_i, \qquad &
    w_{i1} &= -a_i b_i f_i-a_i c_i g_i + b_i g_i + c_i f_i, \qquad &  w_{i2} &= -2_i c_i g_i, \\ 
    w_{i3} &= b_i c_i f_i + c_i^2 g_i - a_i f_i + g_i, \qquad &
    w_{i4} &= -b_i c_i f_i-c_i^2 g_i - a_i f_i + g_i. &
\end{aligned}$}
\end{equation}
Here $a_1,b_1,c_1,d_1,e_1,f_1,g_1,\ldots,a_n,b_n,c_n,d_n,e_n,f_n,g_n$ are
$7n$ real parameters.
\end{proposition}

\begin{proof}
Recall (e.g. from \cite[Example 1.1]{DFRS}) that the Pfaffian of a  skew-symmetric
$4  \times 4$ matrix $(z_{kl})$ is the {\em Pl\"ucker quadric}, which cuts out the
  Grassmannian ${\rm Gr}(2,4)$ in $\PP^5$.
  The threefold $Q$ is the intersection of ${\rm Gr}(2,4)$ with the
 hyperplane $z_{13} + z_{24} = 0$.
  Indeed,  in (\ref{eq:quadricQ2}) we see
  $  z_{13} = p_{i2} $ and  $ z_{24} = -p_{i2}$.
  This identifies our quadric with the  {\em Lagrangian
   Grassmannian}, i.e.~$\, Q \,= \, {\rm LGr}(2,4) \,\subset\, {\rm Gr}(2,4) \,\subset\, \PP^5$.
The points 
${\bf P}_i = (z_{12}: z_{13}: z_{14}: z_{23}: z_{24}:z_{34})$ on $Q$ represent
$2$-dimensional subspaces in $\CC^4$ that are isotropic
with respect to the inner product 
\begin{equation*}
    \Omega \,=\, \begin{bmatrix}
        0 & \rm{Id}_2 \\ -\rm{Id}_2 & 0
    \end{bmatrix}.
\end{equation*}
Symmetric $2 \times 2$ matrices arise as
an affine chart of the Lagrangian Grassmannian $Q$. 
Namely, for $z_{12} = 1$,
 the six $z_{kl}$ are the $2 \times 2$ minors of
 $\begin{footnotesize}
 \begin{bmatrix}
1 & 0 & a & b \\
0 & 1 & b & c 
\end{bmatrix}\!.
\end{footnotesize}$
Note that
$z_{13} + z_{24} = b + (-b) =  0$.

Consider the $\Omega$-isotropic lines in $\PP^3$ that are spanned by the rows of the
$2 \times 4$ matrices
$$
\mathcal{P}_i \,= \,
\begin{bmatrix}
\,1\, & a_i & b_i & c_i \,\, \\
\, d_i\, & e_i & b_i d_i \!+\! c_i e_i & 0 \,\,
 \end{bmatrix}
 \quad {\rm and} \quad
 \mathcal{W}_i \, = \,
 \begin{bmatrix}
\,1 \,& a_i & b_i & c_i\,\, \\
\,f_i\,  & g_i  & b_i f_i \!+\! c_i  g_i & 0 \,\, \\
\end{bmatrix}.
$$
Since $\mathcal{P}_i$ and $\mathcal{W}_i$ have the same first row,
the corresponding points in ${\rm LGr}(2,4)$ span a line in $\PP^5$
that is contained in ${\rm LGr}(2,4)$. Conversely, every line in ${\rm LGr}(2,4)$
arises in this manner. This is because two points in ${\rm Gr}(2,4)$ span a line that is contained in ${\rm Gr}(2,4)$ if and only if the two vector spaces in $\CC^4$ these points represent intersect along a $1$-dimensional subspace. We write the six Pl\"ucker coordinates of
$\mathcal{P}_i$ resp.~$\mathcal{W}_i$  as the entries of a
skew-symmetric $4 \times 4$ matrix $(z_{kl})$ with $z_{13} + z_{24} = 0$.
Each matrix entry $z_{kl}$ is a polynomial in the seven parameters
$a_i,b_i,c_i,d_i,e_i,f_i,g_i$. We now express 
$p_{ij}$ and $w_{ij}$ as linear forms in the $z_{kl}$,
as prescribed in (\ref{eq:quadricQ2}). The result of this transformation
is  the formulas  in (\ref{eq:Knpara}). 
 Our parametrization is actually birational if we view 
 $((a_i,b_i,c_i),(d_i:e_i),(f_i:g_i)) $ as a point
 in $ \CC^3 \times \PP^1 \times \PP^1$.
\end{proof}

We now connect our discussion with the  matrix
varieties $\mathcal{V}_{k,n,r}$ studied in Section \ref{sec2}.

\begin{corollary} \label{cor:para25}
The Gram matrix for the $2n$ vectors in~(\ref{eq:Knpara}) parametrizes
the variety $\mathcal{V}_{2,n,5}$.
\end{corollary}
\begin{proof}
We know from Theorem~\ref{thm:Vmain}
that $\mathcal{V}_{2,n,5}$ is irreducible
and its prime ideal is generated by the $6 \times 6$ minors of $X = (x_{ij})$.
Proposition \ref{prop:abcdefg}
gives us the parametric representation 
\begin{equation}
\label{eq:V2n5para}
  X \, = \,
\begin{footnotesize}
\begin{bmatrix}
p_{10} & p_{11} &  p_{12} &  p_{13} & p_{14} \\
w_{10} & w_{11} & w_{12} & w_{13} & w_{14} \\
p_{20} & p_{21} &  p_{22} &  p_{23} & p_{24} \\
 \vdots & \vdots &  \vdots &  \vdots &  \vdots \\
p_{n0} & p_{n1} &  p_{n2} &  p_{n3} & p_{n4} \\
w_{n0} & w_{n1} & w_{n2} & w_{n3} & w_{n4} \\
\end{bmatrix}
\begin{bmatrix} 
1 & 0 & 0 & 0 & 0 \\
0 & -1 & 0 & 0 & 0 \\
0 & 0 & -1 & 0 & 0 \\
0 & 0 & 0 & 1 & 0 \\
0 & 0 & 0 & 0 & -1 
\end{bmatrix}
\begin{bmatrix} 
p_{10} & w_{10} & p_{20} & \cdots & p_{n0} & w_{n0} \\
p_{11} & w_{11} & p_{21} & \cdots & p_{n1} & w_{n1} \\
p_{12} & w_{12} & p_{22} & \cdots & p_{n2} & w_{n2} \\
p_{13} & w_{13} & p_{23} &  \cdots & p_{n3} & w_{n3} \\
p_{14} & w_{14} & p_{24} & \cdots & p_{n4} & w_{n4} \\
\end{bmatrix}\!.
\end{footnotesize}
\end{equation}
Here  $p_{ij}$ and $w_{ij}$ are the polynomials in (\ref{eq:Knpara}).
The parametrization is a map
$\CC^{7n} \rightarrow \mathcal{V}_{2,n,5}$.
This is ensured by the equations (\ref{eq:PWeqn}).
To see that this map is dominant, we refer to
two equivalent geometric interpretations of (\ref{eq:PWeqn}).
On the one hand, these equations say that
the line spanned by $P_i$ and $W_i$ in $\PP^4$ is contained in ${\rm LGr}(2,4)$.
On the other hand, they say that the
row span of the $2 \times 5$ matrix $\begin{bmatrix} P_i \\ W_i \end{bmatrix}$
is an isotropic plane in $\CC^5$. This row span is a point
in the orthogonal Grassmannian ${\rm OGr}(2,5)$. Since
 the orthogonal group ${\rm O}(5)$ acts transitively on ${\rm OGr}(2,5)$,
the assertion follows from the proof of Theorem \ref{thm:Vmain} for $k=2,r=5$.
\end{proof}

In physics applications, the entries of the Gram matrix $X$
themselves are less relevant. One cares about
  {\em physical quantities}
that are functions of the entries of $X$.
To understand these, one must bear in mind that
the $P$-vectors and the $W$-vectors play different roles.
In conformal field theory, a function that is physical must be
invariant under the translations
\begin{equation}
\label{eq:translations} \quad W_i \,\mapsto \, W_i + \alpha_i P_i, \qquad \alpha_i \in \CC.
\end{equation}
See \cite[equation (3.27)]{CPPR}.
We here work with the {\em fundamental invariants} in (\ref{eq:structures}).
They are referred to as
``the basic building blocks'' after \cite[equation (2.6)]{BMPR}.
Their number is 
\begin{equation}
\label{eq:Nnumber}
 \begin{matrix} N \,= \,n^2(n-1)/2
 \,=\, \binom{n}{2} + \binom{n}{2} + n \binom{n-1}{2} . \end{matrix}  
\end{equation}
The basic building blocks  are the
   rational functions (\ref{eq:fundinv}) on our varieties $\mathcal{V}_{2,n,r}
\subset \CC^M$.

Following \cite{BMPR}, our  primary interest is in conformal correlators.
The simplest conformal correlators
depend on only two of the $n$ fields. According to \cite{CPPR},
 the {\em two-point function}  is
\begin{equation}
\label{eq:2point}
\frac{H_{ij}}{P_{ij}^3} \,\,= \,\,
\frac{(W_i \cdot W_j)(P_i \cdot P_j) - (W_i \cdot P_j) (W_j \cdot P_i)}{(P_i \cdot P_j)^3}
\qquad \hbox{for} \,\,
1 \leq i < j \leq n. 
\end{equation}
We define the {\em conformal two-point variety}
 $\,\mathcal{C}^{(2)}_{n,5}$ to be the closure of the image of the map
 \begin{equation}
 \label{eq:C2n5} \mathcal{V}_{2,n,5} \,\dashrightarrow \, \PP^{\binom{n}{2}-1}.
 \end{equation}
 whose coordinates are given by the formula  (\ref{eq:2point}).
 The following was found computationally:
 
  \begin{conjecture} \label{conj:Kndim}
   The dimension of the two-point variety $\,\mathcal{C}^{(2)}_{n,5}$ equals
 $\, {\rm min}\bigl( 5n-11,\binom{n}{2}-1 \bigr)$. 
 \end{conjecture}

 \begin{proposition} \label{prop:Kndim}
 Conjecture \ref{conj:Kndim} holds for $n \leq 30$.
  This implies that the two-point variety
 $\mathcal{C}^{(2)}_{n,5}$ has positive codimension for $n \geq 9$.
In particular, for $n=9$, it is a hypersurface in $\PP^{35}$.
 \end{proposition}

\begin{proof}[Numerical proof]
The dimension is computed numerically from the
rank of the Jacobian matrix of (\ref{eq:C2n5}) evaluated at random points on $\mathcal{C}_{n,5}^{(2)}$ sampled from a normal distribution and the parameterization from Corollary \ref{cor:para25}.
This computation implies the claim about positive codimension for $9 \leq n \leq 30$.
For $n \geq 31$, it holds because
$\binom{n}{2}-1>7n-10 =  {\rm dim}(\mathcal{V}_{2,n,5})$.
 The dimension equation for $\mathcal{V}_{2,n,5}$ can be seen from Theorem \ref{thm:Vmain}.
 \end{proof}

At present we do not know the degree of the
  hypersurface $\mathcal{C}^{(2)}_{9,5} \subset \PP^{35}$.
It is a considerable challenge to compute its defining polynomial.
Section \ref{sec5} reports on steps in that direction.

Following \cite[Section 6]{RSS},
it would also makes sense to consider the variety
$\mathcal{C}^{(3)}_{n,5}$, which is parametrized by both
two-point functions  and three-point functions. The latter are
defined in \cite[equation (2.8)]{BMPR}.
The three-point variety $\mathcal{C}^{(3)}_{n,5}$  lives
in the tensor space $\PP^{\binom{n}{2}-1} \times
\PP^{n \binom{n-1}{2}-1}$.
This is certainly even more challenging to describe, 
and we leave it for future research.

In the present paper, we focus on the relations among the
basic building blocks~(\ref{eq:fundinv}).
To describe these, we start with the homogeneous coordinate
ring of  $(\PP^4)^n \times (\PP^4)^n$.
This is the polynomial ring generated by the coordinates
of the vectors $P_i$ and $W_i$.
This ring has a natural $\ZZ^{2n}$-grading
given by ${\rm degree}(p_{ij}) = e_i$ and ${\rm degree}(w_{ij}) = f_i$,
where $\{e_1,\ldots,e_n,f_1,\ldots,f_n\}$ denotes the standard basis of $\ZZ^{2n}$.
The expressions in (\ref{eq:fundinv}) and (\ref{eq:structures}) live
in the field of fractions of our polynomial ring, and they
are homogeneous in the $\ZZ^{2n}$-grading.
Namely, we have
 \begin{equation}
 \label{eq:grading} {\rm degree}(P_{ij}) = e_i+e_j \,, \,\,
 {\rm degree}(H_{ij}) = e_i+e_j + f_i+f_j \,, \,\,
 {\rm degree}(V_{i,jk}) = e_i + f_i . 
\end{equation}

Now, the relations among the
basic building blocks~(\ref{eq:fundinv})
live in another polynomial ring, denoted $\CC[P,H,V]$.
That polynomial ring has $N $ variables,
namely $\binom{n}{2}$ variables $P_{ij}$,
$\binom{n}{2}$ variables $H_{ij}$, and
$n \binom{n-1}{2}$ variables $V_{i,jk}$,
with the $\ZZ^{2n}$-grading given in (\ref{eq:grading}).
We are interested in the image of the  rational map
that is defined by (\ref{eq:fundinv}). Geometrically, this is a~map
\begin{equation}
\label{eq:nrmap}
 \mathcal{V}_{2,n,r} \,\dashrightarrow\, \CC^N . 
 \end{equation}
 In the cosmology situation described above, we had $r=5$.
 From the algebra perspective, it makes perfect sense to now replace 
 $(\PP^4)^n \times (\PP^4)^n$ with
 $(\PP^{r-1})^n \times (\PP^{r-1})^n$.
 The map (\ref{eq:nrmap}) is defined for any rank $r$ by
  the rational functions on $\mathcal{V}_{2,n,r}$ given in (\ref{eq:fundinv}).
 The set of polynomials  vanishing on the image of (\ref{eq:nrmap}) is a radical
   ideal $I_{n,r}$ in $\CC[P,H,V]$. Note that $I_{n,r}$
  is homogeneous in the $\ZZ^{2n}$-grading.  
We shall study the ideal $I_{n,r}$
and its variety $V(I_{n,r})$ in~$\CC^N$.
  If $r \geq 5$ then $I_{n,r}$ is prime by Theorem \ref{thm:Vmain}.
We begin with a formula for  the dimension.

\begin{corollary} \label{cor:nonRedGrpAction}
Let $n \geq 3$ and $r \geq 4$. Then  the variety $V(I_{n,r})$  has
dimension $(2r-4) n - \binom{r}{2}$.
In particular, this dimension equals $6n - 10$ for $r=5$, and it equals
$4n-6$ for $r=4$. 
\end{corollary}

\begin{proof}
Theorem \ref{thm:Vmain} tells us that
$\dim(\mathcal{V}_{2,n,r}) = (2r-3)n-\binom{r}{2}$,
for $r \geq 5$. The same holds for $r=4$ by 
Remark \ref{rmk:severalsteps}. See also Proposition \ref{prop:pureof}.
We consider the action of the $n$-dimensional additive group $\mathbb{G}_a^n$ on $\mathcal{V}_{2,n,r}$ that is induced by (\ref{eq:translations}).
Explicitly, this action is as follows:
\[
    X \mapsto U_\alpha \ X \ U_\alpha^T, \quad \text{for } X \in \mathcal{V}_{2,n,r} \text{ and } \alpha = (\alpha_1, \dots, \alpha_n) \in \CC^n,
\]
where $U_\alpha$ is the $2n \times 2n$ block-diagonal matrix:
\[
U_\alpha := {\rm diag}\left(\begin{bmatrix}1 & 0 \\ \alpha_1 & 1 \end{bmatrix}, \begin{bmatrix}1 & 0 \\ \alpha_2 & 1 \end{bmatrix},  \dots, \begin{bmatrix}1 & 0 \\ \alpha_n & 1 \end{bmatrix}\right).
\]
The orbit of a generic point in $\mathcal{V}_{2,n,r}$ under this group action is $n = \dim(\mathbb{G}_a^n)$. So, by Rosenlicht’s theorem, any ring generated by invariants has dimension at most $\dim(\mathcal{V}_{2,n,r}) - n$.
The literature on conformal field theory states that $P,H,V$ 
suffice to express all invariants  \cite[Section~4.2.2]{CPPR}. A precise algebraic statement
is that $P,H,V$ generate the field of invariants $\CC(\mathcal{V}_{2,n,r})^{\mathbb{G}_a^n}$.
This is proved in  \cite[Section 4]{BKMP}.
Using a corollary of \cite[Lemma~2.4]{PV}, we now conclude that
this field has the desired transcendence degree over $\CC$ i.e. $\dim(\mathcal{V}_{2,n,r}) - n$.
This degree is exactly the dimension of $V(I_{n,r})$.
\end{proof}

For the rest of this section, we 
fix  $r=5$, the case  of interest in cosmology.
Here $I_{n,r}$ is a prime ideal.
From (\ref{eq:Knpara}) we get a 
rational parametrization of
the irreducible variety $V(I_{n,5})$.

\begin{remark}
For $n=3$, the ideal $I_{3,5}$ is principal. Its generator is the polynomial in
(\ref{eq:3gen}).
\end{remark}

With this, a big goal is to determine the ideal $I_{4,5}$ of codimension $10$ in $N=24$ variables.

\begin{lemma}\label{le:Vmatrixideal}
The prime ideal $I_{4,5}$ contains the $3 \times 3$ minors of the 
$4 \times 3$ matrix 
\begin{equation}
\label{eq:Vmatrix}
V \,\, = \,\, 
\begin{bmatrix}
    \, V_{1,23} & V_{1,24} & \,\,V_{1,34}\\
    \, V_{2,14} & V_{2,13} &\! -V_{2,34}\\
    \, V_{3,14} & V_{3,24} &\! -V_{3,12}\\
    \, V_{4,23} & V_{4,13} & \,\,V_{4,12}
\end{bmatrix}. 
\end{equation}
The product
$ \,V \cdot
\begin{bmatrix}
    P_{23}P_{14} & - P_{24}P_{13} & P_{34}P_{12}
\end{bmatrix}^T $ is zero on $V(I_{4,5})$, so
its four entries are in $I_{4,5}$.
The eight cubic polynomials above
generate a prime ideal of codimension $4$ and degree~$33$. 
\end{lemma}

\begin{proof}  One can check that the $4 \times 3$ matrix $V$ has rank $2$
under the parametrization of $V(I_{4,5})$ given by (\ref{eq:Knpara}), and that the vector
$\begin{bmatrix}    P_{23}P_{14} & - P_{24}P_{13} & P_{34}P_{12}\end{bmatrix}^T$
 lies in the kernel of~$V$. Using  {\tt Macaulay2} \cite{M2}, we verified
  that the ideal is prime of codimension $4$ and degree $33$.
\end{proof}

We can embark towards our ``big goal'' by extending
 the relation (\ref{eq:3gen})
from $n=3$ to $n=4$. Replacing the indices with any triple
from $\{1,2,3,4\}$, we obtain four minimal generators of
$I_{n,5}$. However, these are far from sufficient to represent
the rank $5$ constraint on the Gram matrix $X$. Indeed,
the $6 \times 6$ minors of $X$ give rise to much more complicated
relations among the quantities in (\ref{eq:fundinv}). We will
study these relations in Section \ref{sec5}. 
In Section \ref{sec4}
we address a lower-dimensional version of the same problem,
namely we consider $d=2$ instead of $d=3$.

\section{Isotropic Vectors in 4-Space}
\label{sec4}

The main result in this section is
Theorem \ref{thm:twopointvariety}, which is
the $d=2$ analogue to Conjecture~ \ref{conj:Kndim}.
Before getting there, we 
study  relations among the
basic building blocks (\ref{eq:fundinv}) for~Gram matrices
of rank $r = d+2 = 4$.  One key difference to $d=3$ is that
the variety $\mathcal{V}_{2,n,4}$ is reducible.
We will restrict our attention to the distinguished
component, which is identified next.

    \begin{proposition} \label{prop:pureof}
    The variety $\mathcal{V}_{2,n,4}$ is pure of dimension $5n-6$.
    It has a main irreducible component $\widetilde{\mathcal{V}}_{2,n,4}$,
    defined by the property that all $\binom{n}{2}$ blocks are nonsingular $2 \times 2$ matrices. 
 \end{proposition}
 
 \begin{proof}
 Purity is Corollary \ref{cor:pure} for $k=2,r=4$.
 The component is $\widetilde{\mathcal{V}}_{2,n,4}=\varphi(\mathcal{Y}_{I,J})$ for
 $I = \emptyset$ and $J = [n]$ in the proof of Theorem \ref{thm:components}.
The $2\times 2$ blocks are nonsingular by  Lemma \ref{lem:XY^T}.
    \end{proof}
 
 As  it was defined in Section \ref{sec3},
 the ideal  $I_{n,4}$  is not prime.
 We now change that definition.
    Namely, we redefine $I_{n,4}$ to be the prime ideal 
of all polynomials that vanish on the image of
 \begin{equation}
\label{eq:n4map}
\widetilde{\mathcal{V}}_{2,n,4} \,\dashrightarrow\, \CC^N . 
 \end{equation}
  We shall study the
$\ZZ^{2n}$-graded prime ideal $I_{n,4}
\subset \CC[P,H,V]$
and its variety $V(I_{n,4}) \subset \CC^N$.
Our first step is to give a parametrization of
$\widetilde{\mathcal{V}}_{2,n,4}$ that works well
in our physics context. 
The formulas in (\ref{eq:uvminors}) below are written in the notation of the
{\em spinor-helicity formalism},
which is widely used in the study of
scattering amplitudes. See \cite{EPS} and the references
therein.

We consider $3n$ unknown vectors
$u_s = (u_{s1},u_{s2})$,
$\overline{u}_s = (\overline{u}_{s1},\overline{u}_{s2})$
and $v_s = (v_{s1},v_{s2})$
in $\CC^2$.
The following notation is used for $2 \times 2 $ minors
of the $3n \times 2$ matrix with these row vectors:
\begin{equation}
\label{eq:uvminors}
 [\,i \,j\, ] = u_{i1} u_{j2} - u_{i2} u_{j1} \,, \,\,
[\,i \, \overline{s} \,] = u_{i1} \overline{u}_{s2}  -  u_{i2} \overline{u}_{s1} \, , \,\,
\langle\, i \,j\, \rangle  =  v_{i1} v_{j2} - v_{i2} u_{j1} .
\end{equation}
Here and in the following proposition, the incides
$i,j,s,t$ are  taken from $[n] = \{1,2,\ldots,n\}$.

\begin{proposition} \label{prop:spinorhelicity}
The variety $\widetilde{\mathcal{V}}_{2,n,4}$ has the following parametric
representation:
\begin{equation}
\label{eq:PHVpara}
P_{ij} \,=\, [\,i\,j\,] \langle \,i\, j\, \rangle \,, \quad
H_{st} \,=\, \langle s\, t \rangle^2  [\,s\,\overline{s}  \,] [\,t\,\overline{t}\, ] \, , \quad
V_{s,ij} \,= \,\frac{\langle s \, i \, \rangle 
 \langle \,s \, j \, \rangle [ \,s \,\overline{s}\, ] 
}{ \langle\, i\, j \,\rangle }  .
\end{equation}
\end{proposition}

\begin{proof}
The following matrix product for $X$
gives a parametrization of the variety
$\widetilde{\mathcal{V}}_{2,n,4} $:
$$
\begin{small} 
\begin{bmatrix}
 u_{11} v_{11} & u_{11} v_{12} & u_{12} v_{12} & \! - u_{12} v_{11} \\
\overline{u}_{11} v_{11} & \overline{u}_{11} v_{12} & 
\overline{u}_{12} v_{12} & \! - \overline{u}_{12} v_{11} \\
   \vdots & \vdots & \vdots & \vdots \\
 u_{n1} v_{n1} & u_{n1} v_{n2} & u_{n2} v_{n2} & \! - u_{n2} v_{n1} \\
\overline{u}_{n1} v_{n1} & \overline{u}_{n1} v_{12} & 
\overline{u}_{n2} v_{n2} & \!- \overline{u}_{n2} v_{n1} \\
\end{bmatrix}
\end{small}
\!
\begin{bmatrix}  
 0 & 0 & 1 & 0 \\   
 0 & 0 & 0 & 1 \\   
 1 & 0 & 0 & 0 \\   
 0 & 1 & 0 & 0 
 \end{bmatrix}
 \!
 \begin{small}
 \begin{bmatrix}
 u_{11} v_{11} &  \overline{u}_{11} v_{11} &    \cdots &
  u_{n1} v_{n1} &  \overline{u}_{n1} v_{n1} \\
 u_{11} v_{12} &  \overline{u}_{11} v_{12} &    \cdots &
  u_{n1} v_{n2} &  \overline{u}_{n1} v_{n2} \\
  u_{12} v_{12} &  \overline{u}_{12} v_{12} &    \cdots &
  u_{n2} v_{n2} &  \overline{u}_{n2} v_{n2} \\
-u_{12} v_{11} & \! \!\!-\overline{u}_{12} v_{11} & \cdots &
 \!\! \!- u_{n2} v_{n1} &\! \!\! -\overline{u}_{n2} v_{n1} \\
 \end{bmatrix}\! .
 \end{small}
 $$
One checks that the $2 \times 2$ diagonal blocks are zero,
and the other $2 \times 2$ blocks are invertible. 
This parametrization of $\widetilde{\mathcal{V}}_{2,n,4} $ is
analogous to that for $\mathcal{V}_{2,n,5}$ given below
Corollary \ref{cor:para25}. Substituting the entries $x_{ij}$ of the  matrix $X$ into 
(\ref{eq:fundinv}), and using (\ref{eq:uvminors}), we obtain (\ref{eq:PHVpara}).
\end{proof}
    
We are now prepared to study the prime ideal $I_{n,4}$.
From the parametrization (\ref{eq:PHVpara}), we~see
$$
V_{s,kt} \, = \,\begin{small}
\frac{\langle s \, k \, \rangle 
 \langle \,s \, t \, \rangle [ \,s \,\overline{s}\, ] 
}{ \langle\, k \, t \,\rangle }  \end{small} \quad {\rm and} \quad
 V_{t,ks} \,= \,\begin{small}
\frac{\langle t \, k \, \rangle 
 \langle \,t \, s \, \rangle [ \,t \,\overline{t}\, ] 
}{ \langle\, k \, s \,\rangle }  . \end{small}
$$
The product of these two expressions equals
$ \,V_{s,kt} V_{t,ks} = -  \langle \,s \, t \, \rangle^2 [ \,s \,\overline{s}\, ] [ \,t \,\overline{t}\, ] 
 =  - H_{st}$.
Therefore,
\begin{equation} \label{eq:HVrelation}
\qquad H_{st} \, + \, V_{s,kt} V_{t,ks} \, \in \, I_{n,4} 
\qquad {\rm for} \,\,
1 \leq s < t \leq n \,\,\,{\rm and}\,\,\,
k \in [n] \backslash \{s,t\} .
\end{equation}
This means that  $I_{n,4}$ is determined
by its elimination ideal $I_{n,4} \cap \CC[P,V]$.
A computation shows that this elimination ideal is zero for $n=3$. This 
brings us back to the Introduction:

\begin{remark} The ideal in (\ref{eq:threebinomials}) is equal to $I_{3,4}$.
Hence $V(I_{3,4})$ is a complete intersection in~$\CC^9$.
\end{remark}

Fix $n \geq 4$.
Our task is to find all  relations among the 
$\binom{n}{2} + n\binom{n-1}{2}$
invariants $P_{ij}$ and $V_{i,jk}$.
The $\binom{n}{2}$ invariants $H_{ij}$ are  determined by (\ref{eq:HVrelation}).
Together, these relations generate $I_{n,4}$.
The next result solves our problem for $n=4$. The first two statements 
mirror Lemma \ref{le:Vmatrixideal}.
 
\begin{proposition}\label{prop:I44}
The ideal $I_{4,4}$ contains the $2 \times 2$ minors of the 
$4 \times 3$ matrix $V$ in~(\ref{eq:Vmatrix}).
The product
$ V \cdot
\begin{bmatrix}
    P_{23}P_{14} & - P_{24}P_{13} & P_{34}P_{12}
\end{bmatrix}^T $ is zero 
for (\ref{eq:PHVpara}), so
its four entries are in $I_{4,4}$.
Furthermore, the ideal $I_{4,4}$ contains
the $10$ distinct entries of the symmetric $4 \times 4$ matrix
\begin{equation} \label{eq:new4by4}
V \cdot\begin{bmatrix}
    0 & 1 & - 1\\
    1 & 0 & 1\\
    -1 & 1 & 0 
\end{bmatrix} \cdot V^T .
\end{equation}
The above relations, together with the six quadrics in (\ref{eq:HVrelation}),
generate the prime ideal $\,I_{4,4}$.
\end{proposition}

\begin{proof}
This is proved by computations with {\tt Macaulay2} \cite{M2}.
\end{proof}

\begin{corollary} \label{coro:In4} Let $n \geq 5$.
For every quadruple $\{i_1,i_2,i_3,i_4\}$ in $[n]$,
by relabeling the polynomials in Proposition \ref{prop:I44},
we obtain a subset of the minimal generators of the ideal $I_{n,4}$.
\end{corollary}

This set of polynomials does not suffice to generate the 
prime ideal $I_{n,4}$. But, we conjecture that $I_{n,4}$ is a
minimal prime of the ideal of relations in Corollary \ref{coro:In4}.
This would follow if we could show that this ideal has dimension
$4n-6$. Indeed, this is  the dimension of $I_{n,4}$,
by Corollary \ref{cor:nonRedGrpAction}.
We verified this computationally for $n \leq 10$.
We now zoom in on $n=5$.

\begin{proposition}
The prime ideal $I_{5,4}$ is minimally generated by
$180$ quadrics, $20$ cubics, $95$ quartics and $156$ quintics
in $50$ variables.
Its variety $V(I_{5,4})$  has dimension $14$ and degree $10145$.
The relations in Corollary \ref{coro:In4} define a reducible variety in $\CC^{50}$ 
with components of various dimensions.  It has dimension $14$, and hence
$V(I_{5,4})$ is a top-dimensional component.
\end{proposition}
   
\begin{proof}[Proof and Discussion]
Let $J$ be the ideal generated  by the $2 \times 2$ minors of
 all $3 \times 4$ matrices obtained by relabeling (\ref{eq:Vmatrix})
 and the entries of all $4 \times 4$ matrices
 obtained by relabeling (\ref{eq:new4by4}).
The ideal $J$ is generated by $130$ quadrics.
It has codimension $23$, equal to the codimension of $I_{5,4} \cap \CC[V]$.
 It is not prime, in fact it has $23$ minimal primes. 
 Every minimal prime, except for one, contains at least one of the
 variables $V_{i,jk}$. By saturating $J$ with respect to all
 variables, we obtain the distinguished associated prime $\tilde J$,
 which is equal to       $I_{5,4} \cap \CC[V]$. 
 The prime ideal $\tilde J$ has codimension $23$, degree $130$ and is 
 minimally generated by $170$ quadrics.

We now consider the ideal in $\CC[P,V]$ which is generated by $\tilde J$
together with all relations  that are obtained by relabeling
$\, V \!\cdot\! [  P_{23}P_{14}, - P_{24}P_{13} , P_{34}P_{12}]^T$.
    That larger ideal is not prime, but it has $42$ minimal primes of various codimensions. We saturate with respect to the variables again, and obtain
    one distinguished minimal prime of codimension $26$ and degree $835$. It is equal to
 $I_{5,4} \,\cap\, \CC[P,V]$ and is generated by $170$ quadrics, $20$ cubics,
 $95$ quartics and $156$ quintics. 
We now add the $10$ quadrics in (\ref{eq:HVrelation}) to obtain
$I_{5,4}$. This ideal is prime, and it has
  codimension $36$ and degree $10145$, as was asserted.
\end{proof}

Even more interesting is the two-point variety 
$\mathcal{C}^{(2)}_{n,4} $ and the three-point variety $\mathcal{C}^{(3)}_{n,4}$.
The former lives in the matrix space $\PP^{\binom{n}{2}-1}$.
The latter lives in the tensor space
$\PP^{\binom{n}{2}-1} \times \PP^{n \binom{n-1}{2}-1}$.
The definitions are as in Section \ref{sec3} but with $d=2$.
Following (\ref{eq:2point}), the two-point functions~are
\begin{equation}\label{eq:paramC2n4}
     z_{st} \,\, := \,\, \frac{H_{st}}{P_{st}^3} \,\, = \,\,
\frac{ [s \, \overline{s}]\, [t \, \overline{t}]
}{ [ s \, t ]^3 \langle s \, t \rangle} 
\qquad {\rm for} \,\,\, 1 \leq s < t \leq n .
\end{equation}
The two-point variety $\mathcal{C}_{n,4}^{(2)}$ is the closure of the image
of the  corresponding map
 \begin{equation}
 \label{eq:C2n4} \widetilde{\mathcal{V}}_{2,n,4} \,\dashrightarrow \, \PP^{\binom{n}{2}-1}.
 \end{equation}
 
Our main result  in this section is the formula for the dimension of this variety.

 \begin{theorem} \label{thm:twopointvariety}
   The dimension of the two-point variety $\,\mathcal{C}_{n,4}^{(2)}$ equals $\,3n - 7\,$ for $\,n \ge 3$.
\end{theorem}
 
\begin{proof}
We compute the dimension using tropical geometry \cite{MS}.
The parametrization \eqref{eq:paramC2n4} expresses $\mathcal{C}_{n,4}^{(2)}$ 
as the Hadamard product of the reciprocal variety of
the Grassmannian  $ {\rm Gr}(2,n)$ given by $\langle s \,t \rangle $ and the image of 
${\rm Gr}(2,2n)$ under the monomial map $y_{st} = [s\,\bar{s}] [t\,\bar{t}] / [s\,t]^3$. 
Points on that image satisfy equations of degree $12$  that are derived from Pl\"ucker relations, like
\begin{equation}
\label{eq:cubedgrassmannian}  \begin{matrix}
\,\,\,\,\quad y_{12}^3 y_{13}^3 y_{24}^3 y_{34}^3 
-y_{12}^3 y_{14}^3 y_{23}^3 y_{34}^3
+y_{13}^3 y_{14}^3 y_{23}^3 y_{24}^3 
\,+\,21 y_{12}^2 y_{13}^2 y_{14}^2 y_{23}^2 y_{24}^2 y_{34}^2 \\
-\,3 y_{12}^3 y_{13}^2 y_{14} y_{23} y_{24}^2 y_{34}^3 
+3 y_{12}^3 y_{13} y_{14}^2 y_{23}^2 y_{24} y_{34}^3 
+3 y_{12}^2 y_{13}^3 y_{14} y_{23} y_{24}^3 y_{34}^2 \\
\,+\,3 y_{12}^2 y_{13} y_{14}^3 y_{23}^3 y_{24} y_{34}^2
+3 y_{12} y_{13}^3 y_{14}^2 y_{23}^2 y_{24}^3 y_{34}
-3 y_{12} y_{13}^2 y_{14}^3 y_{23}^3 y_{24}^2 y_{34}.
\end{matrix}
\end{equation}
By \cite[Theorem 4.3.5]{MS},
the tropicalization of  ${\rm Gr}(2,n)$ is the
space $\Delta_n$ of phylogenetic trees on $n$ leaves, and similarly for
${\rm Gr}(2,2n)$ and $\Delta_{2n}$.
Using \cite[Proposition 5.5.11]{MS}, we conclude
\begin{equation}
\label{eq:tropCn4}
    {\rm trop}\bigl(\mathcal{C}_{n,4}^{(2)}\bigr) \,\,= \,\,
    (-\Delta_n) \,+\,E \cdot \Delta_{2n} .
\end{equation}
Here $E$ is the $\binom{n}{2} \times \binom{2n}{2}$ matrix whose rows  are indexed by 
pairs $\{s,t\} \in \binom{[n]}{2}$, where the nonzero entries in row 
$\{s,t\}$ are $-3$ in column $[s\,t]$ and $1$ in columns $[s\,\bar{s}]$ and $[t\,\bar{t}]$. 

Under the linear map $E: \RR^{\binom{2n}{2}} \rightarrow \RR^{\binom{n}{2}}$, the
 lineality space of $\Delta_{2n}$ maps onto the lineality space of $\Delta_n$.
 We now investigate the polyhedral fan on the right in (\ref{eq:tropCn4}).
We claim that
\begin{equation}\label{eq:E.Delta2n}
    E \cdot \Delta_{2n} \,= - \,  \Delta_{n}.
\end{equation}
We deduce this from the structure of the matrix $E$. 
A column is zero if it is indexed by $[\bar{s}\, \bar{t}]$
or by $[s \,\bar{t}]$ with $s \neq t$.
Each column indexed by $[s \,\bar{s}]$ is in the lineality space of $\Delta_n$. 
The remaining columns of $E$ form the matrix $- 3 \,{\rm Id}_{\binom{n}{2}}$. This shows that $\,    E \cdot \Delta_{2n} =  - 3 \cdot \Delta_{n} = - \Delta_{n}$. 

Substituting (\ref{eq:E.Delta2n}) into 
(\ref{eq:tropCn4}), and bearing in mind the lineality space of $\Delta_n$,
we obtain
$$ \dim \bigl(\,{\rm trop}(\,\mathcal{C}_{n,4}^{(2)}\,) \,\bigr) \,=\,
\dim \bigl( \,(-\Delta_n) + (-\Delta_n)\, \bigr) \, = \,
(2n-4) + (2n-4) - (n-1) \, = \,3n - 7.
$$
The proof is now completed by the dimension part in the Structure Theorem \cite[\S 3.3]{MS}.
\end{proof}

 \begin{cor} \label{cor:degree270} The two point-variety $\mathcal{C}_{n,4}^{(2)}$ has positive codimension for $n \ge 5$. In particular, for $n =5$ it is a hypersurface in $\PP^9$.
 This hypersurface has degree $270$ and reciprocal degree~$90$.
\end{cor}

\begin{proof}[Proof and Discussion]
The first two sentences are immediate from Theorem \ref{thm:twopointvariety}.
By the {\em reciprocal degree} we mean the degree of
hypersurface in $\PP^9$, which is parametrized by
\begin{equation}\label{eq:paramC2n4rec}
     \frac{1}{z_{st}} \,\, := \,\, \frac{P_{st}^3}{H_{st}}\,\, = \,\,
\frac{ [ s \, t ]^3 \langle s \, t \rangle}{ [s \, \overline{s}]\, [t \, \overline{t}]}
\qquad {\rm for} \,\,\, 1 \leq s < t \leq 5 .
\end{equation}
This is more convenient than (\ref{eq:paramC2n4}) because the cube is now in the numerator.
To compute the two degrees $90$ and $270$, we applied the numerical techniques
for hypersurface degrees described in \cite[Section 5]{BHORS}.
Specifically, we replaced  (\ref{eq:C2n4}) by a map with only $9$ parameters that is $4$-to-$1$. 
We then sampled $9$ random affine-linear equations on the image,
and we solved the corresponding equations in the $9$ parameters
using {\tt HomotopyContinuation.jl} \cite{julia}.
These equations were found to have $360 = 4 \cdot 90$
and $1080 = 4 \cdot 270$ solutions respectively.

For $n \geq 6$, it will be possible to determine the degree of $\mathcal{C}_{n,4}^{(2)}$ using
the multiplicities on
${\rm trop}(\mathcal{C}_{n,4}^{(2)}) = (-\Delta_n) + (-\Delta_n)\,$
from the proof above.
But we leave this to a future project.
\end{proof}

\section{Isotropic Vectors in 5-Space}
\label{sec5}

We now return to the case $r=d+2 = 5$ which is relevant for
cosmology \cite{BMPR}. Our goal is to find
 prime ideal $I_{n,5}$ of all relations among the
basic building blocks in (\ref{eq:structures}).
We here present our data for
$n=4$ fields.  Computing $I_{4,5}$ turned out
to be a surprisingly hard problem in eliminating variables.
What follows is a case study in computer algebra
which is of independent interest.
It highlights current numerical tools for implicitization
in action.

Our point of departure for $n=4$ is the list of all
$6 \times 6$ minors of the $8 \times 8$ matrix $X$.
Each minor is the determinant of a $6 \times 6$ submatrix 
obtained by deleting two rows $i$ and $j$
and two columns $k$ and $l$, where $i,j,k,l \in \{1,2,\ldots,8\}$.
This minor is denoted by $X_{ij,kl}$. This is an
element of degree six in the polynomial
ring $\CC[x_{13},x_{14},\ldots,x_{68}]$, which has $24$ unknowns.
The minors are homogeneous with respect to the natural
$\ZZ^8$-grading. For instance,
\[
    \begin{matrix}
    \hbox{$X_{12,12}$ is a sum of \ $40$ monomials of degree $(0,0,2,2,2,2,2,2)$,}\\
    \hbox{$X_{12,23}$ is a sum of \ $92$ monomials of degree $(1,0,1,2,2,2,2,2)$,}\\
    \hbox{$X_{12,34}$ is a sum of $150$ monomials of degree $(1,1,1,1,2,2,2,2)$,}\\
    \hbox{$X_{13,24}$ is a sum of $115$ monomials of degree $(1,1,1,1,2,2,2,2)$.}
    \end{matrix}
\]
A symmetric $8 \times 8$ matrix has $406 $ distinct
$6 \times 6$ minors, but these
are linearly dependent.

\begin{proposition}\label{prop:minorsOfM4}
 The $406$ minors  of size $6 \times 6$ in $X$ satisfy the
$70 = \binom{8}{4}$ linear relations
\begin{equation}
\label{eq:surpriseplucker} \qquad
X_{ij,kl} \, - \, X_{ik,jl} \, + \, X_{il,jk} \,\, = \,\, 0 
\qquad {\rm for} \,\, 1 \leq i < j < k < l \leq 8.
\end{equation}
The ideal of $6 \times 6$ minors has $336$ minimal generators.
These occur in $266 =
\binom{8}{4} \!+\! 8 \binom{7}{2} \!+\! \binom{8}{2}$
distinct $\ZZ^8$-degrees. It is the prime ideal of
 $\,\mathcal{V}_{2,4,5}$, which has
 codimension $6$ and degree~$672$.
\end{proposition}

\begin{proof}
The first three sentences
were verified by direct computation.
The relations  (\ref{eq:surpriseplucker})
look like the Pl\"ucker relations for the Grassmannian
${\rm Gr}(2,8)$. An analogous statement holds for all
symmetric matrices and all Grassmannians. We shall explain this in
 Theorem \ref{thm:PluckerRelsForSymmMatrices} below.
The last sentence in Proposition \ref{prop:minorsOfM4} follows from
Theorem \ref{thm:Vmain} and Corollary  \ref{cor:expecteddegree}.
\end{proof}

The following general fact about
linear relations among the $k \times k $ minors of
a symmetric matrix with unknown entries is of independent interest for 
linear algebra. It explains (\ref{eq:surpriseplucker}).

\begin{theorem}\label{thm:PluckerRelsForSymmMatrices}
    Let $X$ be an $m \times m$ symmetric matrix with indeterminate entries, $k \leq m$, and let $I = \{ i_1< \dots < i_{k-1}\}$, $J = \{j_1 < \dots < j_{k+1}\}$ be two subsets of $[n]$ of sizes $k-1$ and $k+1$ respectively. Then the following Pl\"ucker-like relation holds among the minors of $X$:
    \[        \sum_{\ell=1}^{k+1}\, (-1)^{\ell} \ X_{Ij_\ell, J\setminus j_\ell} = 0.    \]
\end{theorem}

\begin{proof}
This identity appears 
in \cite[Lemma 5.2, p.~345]{DCP},
where it is written as $ \dot T = 0$.
Conca \cite[page 414]{Conca} uses the identity 
when he proves that symmetric minors are a Gr\"obner basis.
\end{proof}

Proposition \ref{prop:minorsOfM4} characterizes the prime  ideal of $\mathcal{V}_{2,4,5}$.
From that ideal in $\CC[X]$, we now need to compute the ideal
$I_{4,5}$ in $\CC[P,H,V]$. Note that both polynomial rings have $24$ variables.
This is a problem of computer algebra, which can in principle be solved as follows.
Augment the $336$ minors in Proposition \ref{prop:minorsOfM4} with the
$24$ relations (\ref{eq:fundinv}), and eliminate the unknowns in $X$.
This involves the clearing of the denominators $x_{2j-1,2k-1}$ for $V_{i,jk}$.
The result of this elimination is the desired prime ideal $I_{4,5}$ in the polynomial ring
in $\CC[P,H,V]$.

This sounds nice and easy in theory, but it is very difficult in practice.
We used the multigrading method of Cummings and Hollering \cite{CH},
adapted from polynomials to rational functions.
We  hunt for polynomials in $I_{4,5}$ degree by degree, beginning with the four minors
\begin{equation}
\label{eq:fourminors}
X_{12,12},\,\,
X_{34,34},\,\,
X_{56,56},\,\,
X_{78,78}.
\end{equation}
These are clearly polynomials in the variables $P_{ij}, H_{ij}, V_{i,jk}$.
For instance, $X_{78,78}$
equals (\ref{eq:3gen}).

We first attempted to express the other
 $332$ minors of $X$ as polynomials in the $24$ variables $P,H,V$.
 But this was unsuccessful in all cases.
It is not possible to rewrite any of the other
minors without introducing denominators.
This is a consequence of the following lemma:

\begin{lemma} \label{lem:six}
The ideal $I_{4,5}$ has only $12$ minimal generators in degrees $\leq 6$, namely
the eight cubics in Lemma \ref{le:Vmatrixideal}
and the four sextics in (\ref{eq:fourminors}), obtained by 
relabelling the sextic in (\ref{eq:3gen}).
\end{lemma}

 This was proved by linear algebra as in \cite{CH}.
Lemma \ref{lem:six} refers  to the $\ZZ$-grading  given by
$$ {\rm degree}(P_{ij}) = 1, \,\,{\rm degree}(H_{ij}) = 2, \,\,{\rm degree}(V_{i,jk}) = 1. $$
As in \cite{CH}, we find
 all polynomials in $I_{4,5}$ with a fixed $\ZZ^8$-degree by solving  linear equations. 
 Our ansatz is a linear combination with unknown coefficients of
   all monomials in that degree. We then
   substitute the parametrization from  (\ref{eq:Knpara}) and (\ref{eq:structures})
into the ansatz. Each numerical evaluation of the 
$28$ parameters $a_1,a_2,\ldots,g_4$ furnishes one linear equation
in the unknowns. We form a system of sufficiently many such linear equations.
By solving this linear system, we
obtain a basis for the $\ZZ^8$-graded component of our ideal $I_{4,5}$.
We prove Lemma \ref{lem:six} by performing this computation
for all $\ZZ^8$-degrees whose $\ZZ$-degree is at most six.
The software {\tt Macaulay2} is used to extract minimal ideal
generators from these various vector space bases. 

The next two lemmas were found by techniques that are based on what is described~above.

\begin{lemma} \label{lem:seven}
The ideal $I_{4,5}$ has precisely $172$ minimal generators of total degree seven, $156$ of which are polynomials in $H$ and $V$.
The generators are characterized in the following~table:
$$ \begin{small}
\begin{matrix}
    {\rm degree} & {\rm \# \,\, orbits} & {\rm \# \,\, monomials}
     & {\rm dim}_{\rm gr}({I_{4,5}}) & {\rm \# \,\, mingens} \\
    (2,2,2,1,2,2,2,1) & 4 & 1812 & 285 & 15 \\
    (3,2,1,1,3,2,1,1) & 12 & 1375 & 182 & 5\\
    (3,2,2,0,3,2,2,0) & 12 & 784 & 84 & 2\\
    (3,3,1,0,3,3,1,0) & 12 & 600 & 47 & 1\\
    (3,2,2,2,2,1,1,1) & 4 & 945 & 373 & 1\\
    (3,3,2,1,2,2,1,0) & 12 & 540 & 205 & 1
\end{matrix}
 \end{small} $$
\end{lemma}

Column \# 3 gives the number of monomials in the given degree,
column \# 4 reports the dimension of $I_{4,5}$ in that degree,
and column \#5 measures the subspace of ideal generators.

\begin{example}
There are $15$ generators of $I_{4,5}$ in $\ZZ^8$-degree $(2,1,2,2, 2,1,2,2)$.
Each of them has at least $55$ terms.
The smallest polynomial in that basis, with $55$ terms, is equal to
$$ \begin{scriptsize}
\begin{matrix}
   V_{123} V_{134} V_{213} V_{314}^2 V_{413} V_{412}
 + V_{123} V_{134} V_{234} V_{314}^2 V_{413}^2
 + V_{134}^2 V_{213} V_{314}^2 V_{413} V_{412}
 - V_{134}^2 V_{214} V_{314} V_{324} V_{413} V_{412} \\
 + V_{134}^2 V_{234} V_{314}^2 V_{413}^2
 - V_{134}^2 V_{234} V_{314} V_{324} V_{423} V_{413}
 - 2 H_{13} V_{123} V_{213} V_{314} V_{412}^2
 + 2 H_{13} V_{123} V_{213} V_{314} V_{413} V_{412} \\
 - 2 H_{13} V_{123} V_{234} V_{314} V_{413} V_{412}
 + 2 H_{13} V_{123} V_{234} V_{314} V_{413}^2
 + 2 H_{13} V_{124} V_{214} V_{314} V_{412}^2
 - H_{13} V_{124} V_{214} V_{314} V_{413} V_{412} \\
 + 2 H_{13} V_{124} V_{234} V_{314} V_{423} V_{412}
 - H_{13} V_{124} V_{234} V_{314} V_{423} V_{413}
 + 2 H_{13} V_{134} V_{213} V_{314} V_{413} V_{412}
 - H_{13} V_{134} V_{214} V_{314} V_{413} V_{412} \\
 - H_{13} V_{134} V_{214} V_{324} V_{413} V_{412}
 + 2 H_{13} V_{134} V_{234} V_{314} V_{413}^2
 - H_{13} V_{134} V_{234} V_{314} V_{423} V_{413}
 - H_{13} V_{134} V_{234} V_{324} V_{423} V_{413} \\
 - H_{14} V_{123} V_{213} V_{314}^2 V_{412}
 - H_{14} V_{123} V_{234} V_{314}^2 V_{413}
 - H_{14} V_{134} V_{213} V_{314}^2 V_{412}
 + H_{14} V_{134} V_{214} V_{314} V_{324} V_{412} \\
 - H_{14} V_{134} V_{234} V_{314}^2 V_{413} 
 + H_{14} V_{134} V_{234} V_{314} V_{324} V_{423}
 + H_{34} V_{123} V_{134} V_{213} V_{314} V_{412}
 + H_{34} V_{123} V_{134} V_{234} V_{314} V_{413}\\
 + H_{34} V_{134}^2 V_{213} V_{314} V_{412}
 - H_{34} V_{134}^2 V_{214} V_{324} V_{412} 
 + H_{34} V_{134}^2 V_{234} V_{314} V_{413}
 - H_{34} V_{134}^2 V_{234} V_{324} V_{423}
 + H_{13}^2 V_{213} V_{413} V_{412} \\
 - 2 H_{13}^2 V_{213} V_{423} V_{412}
 + H_{13}^2 V_{213} V_{423} V_{413} 
 + H_{13}^2 V_{214} V_{413} V_{412} 
 - H_{13}^2 V_{214} V_{413}^2 
 + H_{13}^2 V_{234} V_{413}^2 
 - H_{13}^2 V_{234} V_{423} V_{413} \\
 - H_{13} H_{14} V_{213} V_{314} V_{412}
 - H_{13} H_{14} V_{213} V_{314} V_{423} 
 - H_{13} H_{14} V_{214} V_{314} V_{412} 
 + H_{13} H_{14} V_{214} V_{314} V_{413} \\
 + 2 H_{13} H_{14} V_{214} V_{324} V_{412}
 - H_{13} H_{14} V_{234} V_{314} V_{413}
 - H_{13} H_{14} V_{234} V_{314} V_{423}
 + 2 H_{13} H_{14} V_{234} V_{324} V_{423} \\
 + H_{13} H_{34} V_{123} V_{213} V_{412}
 + H_{13} H_{34} V_{123} V_{234} V_{413}
 + H_{13} H_{34} V_{124} V_{214} V_{412} 
 + H_{13} H_{34} V_{124} V_{234} V_{423} \\
 + H_{13} H_{34} V_{134} V_{213} V_{412}
 - H_{13} H_{34} V_{134} V_{214} V_{412}
 + H_{13} H_{34} V_{134} V_{234} V_{413}
 - H_{13} H_{34} V_{134} V_{234} V_{423}.
 \end{matrix}
\end{scriptsize}
$$
Nonlinear algebra is not as simple as the relations in (\ref{eq:3gen}) and
(\ref{eq:threebinomials}) might have suggested.
\end{example}

\begin{lemma} \label{lem:eight}
The ideal $I_{4,5}$ has precisely $796$ minimal generators
of total degree eight; $592$ of them are polynomials in $H$ and $V$.
The generators are characterized in the following table:
$$ \begin{footnotesize}
\begin{matrix}
    {\rm degree} & {\rm \# \,\, orbits} & {\rm \# \,\, monomials} & {\rm dim}_{\rm gr}({I_{4,5}}) & {\rm \# \,\, mingens}\\
    (2,2,2,2,2,2,2,2) &  1 & 4539 & 1048 & 34\\
    (3,2,2,1,3,2,2,1) & 12 & 3426 &  670 & 14\\
    (3,3,1,1,3,3,1,1) &  6 & 2596 &  423 & 10\\
    (3,3,2,0,3,3,2,0) & 12 & 1480 &  211 & 9\\
    (4,2,1,1,4,2,1,1) & 12 & 2209 &  327 & 6\\
    (4,2,2,0,4,2,2,0) & 12 & 1261 &  157 & 6\\
    (4,3,1,0,4,3,1,0) & 24 &  965 &   87 & 3\\
    (4,4,0,0,4,4,0,0) &  6 &  371 &    1 & 1\\
    (3,3,2,2,2,2,1,1) &  6 & 2232 & 1004 & 12\\
    (3,3,3,1,2,2,2,0) &  4 & 1272 &  558 & 6\\
    (4,2,2,2,3,1,1,1) &  4 & 1704 &  732 & 3\\
    (4,3,2,1,3,2,1,0) & 24 &  975 &  406 & 3\\
    (4,4,1,1,3,3,0,0) &  6 &  438 &  163 & 1\\
    (4,3,3,2,2,1,1,0) & 12 &  474 &  209 & 1\\
    (4,4,2,2,2,2,0,0) &  6 &  276 &  117 & 1
\end{matrix}
 \end{footnotesize}
$$
\end{lemma}

By combining Lemma \ref{lem:six}, \ref{lem:seven} and \ref{lem:eight},
we currently have a list of $980$ known minimal generators for
the prime ideal $I_{4,5}$. We do not know yet whether that list is complete.

\begin{proof}[Proof and discussion for Lemmas \ref{lem:seven} and \ref{lem:eight}]
The linear systems in each degree are large. Since it is challenging
to solve them over $\QQ$, we worked over a finite field $\ZZ/p\ZZ$ instead.
Suppose that  $I_{4,5}$ has generators with integer coefficients less
than $p$ in absolute value.
We can then lift any generator $f$ over $\ZZ/p\ZZ$ to a generator over $\QQ$ by finding 
$m \in \{-p+1,\ldots,p-1\}$ such that $m \cdot f$ vanishes on 
 $V(I_{4,5})$. We performed these  linear algebra computations
 for each degree in \texttt{Macaulay2}. In order to succeed,
we experimented with various prime numbers $p$.

Our computations rest on the following combinatorial
considerations. Fix $\ZZ$-degree $7$. We need to explore all
$\ZZ^8$-degrees. We start with  $(\lambda,\lambda) \in \ZZ^8$,
where $\lambda $ is a partition of
$7$ into at most $4$ parts. There are $p_4(7) = 11$ such partitions $\lambda$.
We found that
only the four partitions $(2,2,2,1), (3,2,1,1), (3,2,2,0), (3,3,1,0)$ yield
minimal generators for Lemma \ref{lem:seven}.

We now have all generators in the  variables $H$ and $V$. 
     Next we hunt for irreducible polynomials in $I_{4,5}$ that also contain variables $P$.
     The $\ZZ^8$-degree in $P$ of such a polynomial~is
    $$
    (1,1,1,1,0,0,0,0), \,\,(2,2,2,2,0,0,0,0)\,\text{ or }\,(3,3,3,3,0,0,0,0).
    $$ 
    In the first case, we have degree $2$ 
     in $P$ and degree $5$ in  $(H,V)$.
      There are $6$ partitions of $5$ into at most $4$ parts. Only $2$ of them, namely $(2,1,1,1)$ and $(2,2,1,1)$, contribute minimal generators. We repeat the same for polynomials of degree $4$ and $6$ in the $P$ variables, but we find no minimal generators.
      The same strategy was carried out to establish Lemma \ref{lem:eight}.
\end{proof}

\section{Low Rank Schemes} \label{sec6}

We now depart from physics and return to Gram matrices.
 Section \ref{sec2} was about the
 high rank regime $r > 2k$ where the ideal of  minors is prime.
 In this section, we assume $k \geq 2$, $n \geq 3$ and $0 \leq r \leq kn$.
 Following (\ref{eq:k2r}), we always set
 $ \ell = \lfloor r/2 \rfloor $. Here is our first main result:
 
\begin{theorem} \label{thm:dimension}
The dimension of the variety $\mathcal{V}_{k,n,r}$ of low rank Gram matrices is equal to
\begin{equation}
  \label{eq:truedimension}
  {\rm dim}( \mathcal{V}_{k,n,r}) \,\,\, = \,\,\,
\begin{cases}
   \,\,  n\left(\,\ell (k-\ell) + \ell r - \binom{\ell+1}{2}\, \right) - \binom{r}{2} & {\rm if} \quad r \leq 2k+2,
\smallskip \\   
    \,\,  n\left(\,\,kr - \binom{k+1}{2}\,\,\right) - \binom{r}{2} & {\rm if} \quad r \geq 2k-1.  
    \end{cases}
\end{equation}
In the second case, the ideal of $\,(r+1) \times (r+1)$ minors is Cohen-Macaulay
of expected degree (\ref{eq:expecteddegree}).
Note that the two expressions in (\ref{eq:truedimension}) are equal for $\,r
\in \{2k-1,2k,2k+1,2k+2\}$.
\end{theorem}

\begin{corollary} \label{cor:pure}
If $r \in \{2k-1,2k\}$ then every component of $\,\mathcal{V}_{k,n,r}$
has the same dimension.
\end{corollary}

 The range $8 \leq kn \leq 10$ is summarized in
Table \ref{tab:varieties1}.
 We display the
dimension, degree and number of minimal
generators for the ideal of $(r+1) \times (r+1)$ minors
which defines $\mathcal{V}_{k,n,r}$.
For each entry in Table \ref{tab:varieties1}, the first number
 equals (\ref{eq:truedimension}),
the second is (\ref{eq:expecteddegree}) for $r \geq 2k-1$,
and the third can be derived from Theorem \ref{thm:PluckerRelsForSymmMatrices}.
Recall the number {\bf 336} in
 Proposition~\ref{prop:minorsOfM4}.

\begin{center}
\begin{table}[h]
\vspace{-0.4cm}
$$ \begin{matrix}
\qquad r \,\backslash\, (k,n)\! \!\!\!\!\! \!\!& (2,4) & (2,5) & (3,3) & (4,2) & (5,2) \\
1 & 0, 41, 284 & 0, 121,740 &   0, 28,378 & 0,17,136 & 0 , 26, 325  \\ 
2 & 7, 84, 1092 & 9,600,4770 & 8, 63,1950 & 7,20,416 & 9, 70, 1700  \\
3 & 9, {\bf 4224}, 1764 & 12, 183040, 13860 & 9, 3915,4599 & 7,200,626 & 9,1190,4550 \\
4 &  \! 14, {\bf 2772},1176 & 19, 306735,  19404 & 15, 930, 4977 & 12,20,416 & 16,175,6202\\
5 & 18, 672, {\bf 336} & 25, 151008, 13860 & 17, 9504, 2520 & 12,100,136 & 16, 1750 ,4550\\
6 & 21, 84, 36 & 30, 28314,4950 & 21, 1386, 540 & 15,4,16 & 21, 50, 1700 \\ 
7 & 23, 8, 1 & 39, 2640, 825 & 24, 120, 45 & 15,8,1 & 21 , 250, 325
\end{matrix}
\vspace{-0.5cm}
$$
\caption{The determinantal varieties $\mathcal{V}_{k,n,r}$ for matrix sizes
$kn$ between $8$ and $10$.
\label{tab:varieties1} }
\end{table}    
\end{center}

\vspace{-0.8cm}

 The second theorem in this section concerns the
 irreducible components of our variety.

\begin{theorem}\label{thm:components}
The variety $\,\mathcal{V}_{k,n,r}$ is irreducible if and only if $\,r$ is odd or $r > 2k$.
For~even $4 \leq r \leq 2k$, the
number of irreducible components
of $\,\mathcal{V}_{k,n,r}$ is $\,2^{n-1}$. For $r=2$, it is $\,2^{n-1}-1$.
\end{theorem}

We now embark towards proving the two theorems.
The idea is to construct certain  maps into
$\mathcal{V}_{k,n,r}$. These maps will reveal the dimension
and component structure.
As in Section \ref{sec2}, the local building block is
the ring $\CC[Y]$ of polynomials in the entries of a
$k \times r$ matrix~$Y$.

\begin{proof}[Proof of Theorem \ref{thm:dimension}]
    We begin by proving that $\mathcal{V}_{k,n,2k}$
    is pure. In fact, we argue that the ideal of minors is Cohen-Macaulay.
	We know from Remark \ref{rmk:severalsteps} that the ring $T$ is a complete intersection of
	dimension $k(3k-1)/2$. 
	The ring $\CC[Y]/\langle Y_i Y_i^T : i=1,\ldots,n \rangle $
	is the $n$-fold tensor product of $T$ and hence has
	Krull dimension $nk(3k-1)/2$. 
	The $k(2k-1)$-dimensional Lorentz group
	${\rm SO}(2k)$ acts faithfully on its spectrum,
	and the quotient is the variety $\mathcal{V}_{k,n,2k}$.
	Using (\ref{eq:Mvariables}),
	this argument implies that the codimension of $\mathcal{V}_{k,n,2k} $ 
	in $\CC^M$ is equal to
	\[
	\binom{kn+1}{2}\,-\, n\binom{k+1}{2} \,-\,  n\Big( \,\frac{nk(3k-1)}{2}  
	- k(2k-1)\Big)\,\,=\,\,  \binom{kn-2k+1}{2}.
	\]
	This matches the codimension in  \cite[Proposition~12]{HT} for
	the ideal of $(2k+1) \times (2k+1)$ minors of a generic 
	symmetric $kn \times kn$ matrix. The latter is Cohen-Macaulay.
	Hence so is our linear section with 
		diagonal $k \times k$ blocks equal to zero.
	We conclude that $\mathcal{V}_{k,n,2k}$ is~pure. 
	
	We have now proved the case $r = 2k$ in \eqref{eq:truedimension}.
	The case $r \geq 2k+1$ already appeared in Theorem \ref{thm:Vmain}.
    For the remaining case $r < 2k$, we shall first need the following lemma.
   
    \begin{lemma}\label{lem:Okr}
      The variety $\,\mathcal{T}_{k,r}$ of $\,k\times r$ matrices $Y$ such that $\,YY^T = 0$ has dimension
      \[
        \dim(\mathcal{T}_{k,r })\,\, = \,\,\ell(k-\ell) + \ell r - \binom{\ell+1}{2} \qquad \text{provided }\,\, r \leq 2k
    \,\,\, \text{and}  \,\,\, \,\ell = \lfloor r /2 \rfloor.
      \]
      For $r= 1,3, \dots, 2k-1$, this variety is irreducible but the ideal 
         $\langle Y Y^T\rangle$ is not radical in $\CC[Y]$.
          For $r = 2, 4, \dots, 2k$, it has two irreducible components, 
          and the ideal $\langle Y Y^T\rangle$ is radical.
              \end{lemma}
    
    \begin{proof}
        Fix $r \leq 2k$ and $\ell = \lfloor r/2 \rfloor$. Consider the birational map $\mathcal{T}_{\ell, r} \times \CC^{(k-\ell) \times \ell} \to \mathcal{T}_{k,r}$ which takes  a rank $\ell$ matrix $Z \in   \mathcal{T}_{\ell, r}$ and a $(k-\ell) \times \ell$ matrix $A \in \CC^{(k - \ell) \times \ell}$ to the $k \times r$ matrix obtained by stacking $Z$ on top of $AZ$. From Remark \eqref{rmk:severalsteps} we know that the defining ideal of $\mathcal{T}_{\ell,r}$ is a complete intersection and hence $\dim(\mathcal{T}_{\ell, r}) = \ell r - \binom{\ell + 1}{2}$. Since $\dim(\mathcal{T}_{k,r}) = \dim(\mathcal{T}_{\ell,r}) + \ell(k - \ell)$, we deduce the desired dimension formula.
        If $r = 2\ell$ is even then $\mathcal{T}_{\ell, r} = \mathcal{T}_{\ell, 2\ell}$ decomposes into two  components of the same dimension, so the same holds for $\mathcal{T}_{k,r}$.
        If $r = 2\ell+1$ is odd then $\mathcal{T}_{\ell,r} = \mathcal{T}_{\ell, 2\ell+1}$ is irreducible, 
        hence so is the image $\mathcal{T}_{k,r}$ under our birational map.
    \end{proof}
        
        We now continue our proof of Theorem \ref{thm:dimension} and compute the dimension of $\mathcal{V}_{k,n,r}$ when $r < 2k$. To do so, we proceed as in the proof of Lemma \ref{lem:Okr} and exhibit a birational map 
        \begin{equation}\label{eq:biratMapell}
         \mathcal{V}_{\ell , n, r} \times \left(\CC^{ (k-\ell) \times \ell} \right)^{n} 
         \, \dashrightarrow \,\mathcal{V}_{k,n,r}.
        \end{equation}
      Its input is an $\ell n \times \ell n$ matrix $X \in \mathcal{V}_{\ell, n, r}$ 
      with $\ell \times \ell$ blocks $(X_{ij})_{1 \leq i,j \leq \ell}$ and $n$
       matrices $A_1, \dots, A_n$ of size $(k-\ell) \times \ell$.
       Its output is the
        $kn \times kn$ symmetric matrix with $k \times k$ blocks
        \begin{equation} \label{eq:smartblocks}
        \begin{pmatrix}
            X_{ij} & X_{ij} A_j^T\\
            A_i X_{ij} & A_i X_{ij} A_j^T
        \end{pmatrix}.
\end{equation}
      We deduce that $\dim(\mathcal{V}_{k,n,r}) = \dim(\mathcal{V}_{\ell,n,r}) + n \ell(k-\ell)$. Since $\ell = \lfloor r/2 \rfloor$, we have $2\ell \leq r$. Hence we can apply \eqref{eq:truedimension} to $\mathcal{V}_{\ell, n, r}$ and deduce that  $\dim(\mathcal{V}_{k,n,r})  = n \left( \ell(k-\ell) + \ell r - \binom{\ell+1}{2} \right) - \binom{r}{2}$. Finally, a calculation 
      with binomial coefficients checks that \eqref{eq:truedimension} holds in all cases.
\end{proof}

To establish the component structure in Theorem \ref{thm:components}, we start with the following lemma.

\begin{lemma}\label{lem:XY^T}
    Let  $X$ and $Y$ be generic matrices in $\mathcal{T}_{k,2k}$. The
    $k \times k$ matrix $X Y^T$ is  singular when $k$ is odd and $X,Y$ are in the same component of $\mathcal{T}_{k,2k}$ or when $k$ is even and $X,Y$ are in different components of  $\mathcal{T}_{k,2k}$.
    Otherwise, $X Y^T$ is invertible. Here is a summary:
        \begin{table}[H]
        \centering
        \begin{tabular}{ c c c }
        
          & $X,Y$ in same component  &  $X,Y$ in different components \\
        
       $k$ odd   & Singular  & Generically invertible  \\
        
       $k$ even  & Generically invertible  & Singular  \\
        
        \end{tabular}
        \caption{Behavior of the product $XY^T$ for two matrices $X,Y \in \mathcal{T}_{k,2k}$.}
        \label{tab:my_label}
    \end{table}
\end{lemma}
\begin{proof}
    As in Remark \ref{rmk:severalsteps}, we use the quadratic form
     given by the $2k \times 2k$ diagonal matrix $Q = {\rm diag}(1,-1,\dots, 1,-1)$. Recall from \eqref{eq:twocomponents} that $\mathcal{T}_{k,2k}$ has the two irreducible components
         \begin{equation}
         \label{eq:Ocomponents} \begin{matrix}
        \mathcal{T}^{\pm}_{k,2k} \,\,= \,\,\bigl\{ Y \, \in \mathcal{T}_{k,2k} \,:\,  \det(Y_A) 
        \,=\, \pm \det(Y_{A^c})\,\,\, \text{for all } A \in \binom{[2k]}{k}\bigr\}.
    \end{matrix} 
    \end{equation}
    Now fix $X,Y \in \mathcal{T}_{k,2k}$.
    The Cauchy-Binet Theorem yields
$$        \det(XQY^T)\,\, = \sum_{A \in \binom{[2k]}{k}}  \det(X_{A}) \det((YQ)_{A})  
                  \,\,  = \sum_{A \in \binom{[2k]}{k}}  (-1)^{|A \cap 2\ZZ|} \det(X_{A}) \det(Y_{A}).
$$                    
    If $X$ and $Y$ come from the same component in (\ref{eq:Ocomponents}) then
    \begin{align*}
     \det(XQY^T) \,\,&= \sum_{A \in \binom{[2k]}{k}} \! (-1)^{|A^c \cap 2\ZZ|} \det(X_{A}) \det(Y_{A}) \\
                 &= \,\,(-1)^k \!\!\sum_{A \in \binom{[2k]}{k}}  (-1)^{|A \cap 2\ZZ|} \det(X_{A}) \det(Y_{A}) 
                 \,\,\,= \,\,\,(-1)^k\det(XQY^T).
    \end{align*}
    Similarly, if $X$ and $Y$ belong to different components then
    \[
        \det(XQY^T) \,\,=\,\, (-1)^{k+1} \det(XQY^T).
    \]
This explains the alternating behavior seen in Table \ref{tab:my_label},
and it completes the proof.
\end{proof}

\begin{proof}[Proof of Theorem \ref{thm:components}]
     We use the quadratic form $Q$ from the proof of Lemma \ref{lem:XY^T}. 
     The variety $\mathcal{V}_{k,n,r}$ is irreducible when $r > 2k$, by Theorem \ref{thm:Vmain}.
     Let $r \leq 2k$ be odd. We know from Lemma \ref{lem:Okr} that the variety $\mathcal{T}_{k,r}$ is irreducible. Then $\mathcal{V}_{k,n,r}$ is also irreducible because it is the image of the $n$-fold product $\mathcal{T}_{k,r} \times \dots \times \mathcal{T}_{k,r}$ under the map $(Y_1,\dots, Y_n) \mapsto (Y_iQY_j^T)_{1 \leq i,j \leq n}$.

     Now assume $4 \leq r = 2\ell \leq 2k$. Let $\varphi$ denote the map from $\mathcal{T}_{\ell,r} \times \dots \times \mathcal{T}_{\ell, r}$ onto $\mathcal{V}_{\ell,n,r}$ that maps $(Y_1, \dots, Y_n)$ to the $\ell n  \times \ell n$ symmetric matrix $X = (Y_iQY_j^T)_{1 \leq i,j \leq n}$.
    For any partition $I \sqcup J$ of $[n]$, write $\mathcal{Y}_{I,J}$ for the irreducible component of $\mathcal{T}_{\ell,r} \times \dots \times \mathcal{T}_{\ell,r}$ where the matrices indexed by $I$ 
    are in the component $\mathcal{T}^{+}_{\ell,r}$ and those indexed by $J$ 
    are in the component $\mathcal{T}^{-}_{\ell, r}$.
    We have $\varphi(\mathcal{Y}_{I,J}) = \varphi(\mathcal{Y}_{J,I})$  because there exists an involution $P \in {\mathrm O}(r)$ which switches the two components  $\mathcal{T}^{\pm}_{\ell, r}$.
    Therefore, if $(Y_1, \dots, Y_n) \in \mathcal{Y}_{I,J}$ then $(Y_1P, \dots, Y_nP) \in \mathcal{Y}_{J,I}$,
     and 
     \[    \varphi(Y_1P, \dots, Y_nP) \,\,= \,\,\Big(Y_{i} P \ Q \ (Y_{j}P)^T\Big)_{1 \leq i,j \leq n} \! = 
         \,\,\,\Big(Y_{i} Q Y_{j}^T\Big)_{1 \leq i,j \leq n} \! = \,\,\,\varphi(Y_1, \dots, Y_n).
     \]

The $2^{n-1}$ varieties $\varphi(\mathcal{Y}_{I,J})= \varphi(\mathcal{Y}_{J,I})$
are distinct by  Lemma \ref{lem:XY^T}. We claim that they
are the irreducible components of $\mathcal{V}_{\ell,n,r}$.
For this, it suffices to show that each of them has dimension
$n(\ell r - \binom{\ell+1}{2}) - \binom{r}{2}$, which is the number in
\eqref{eq:truedimension}. This holds because
$\varphi(\mathcal{Y}_{I,J})$ is the quotient of $\mathcal{Y}_{I,J}$ by the group ${\rm SO}(r)$ and the action of this group on $\mathcal{Y}_{I,J}$ is faithful.

To pass from $\ell$ to $k$, we now use 
       the map $\mathcal{V}_{\ell , n, r} \times \left(\CC^{ (k-\ell) \times \ell} \right)^{n}  \dashrightarrow \mathcal{V}_{k,n,r}$ from \eqref{eq:biratMapell} and
\eqref{eq:smartblocks}. Since this map is birational, we conclude that
          $\mathcal{V}_{k,n,r}$ has $2^{n-1}$ irreducible components.         
     Finally, for $r = 2$, recall from Proposition \ref{prop:k=1decomp} that $\mathcal{V}_{1,n,2}$ has $2^{n-1}-1$ components. The birational map \eqref{eq:biratMapell} shows that $\mathcal{V}_{k,n,2}$ 
      has $2^{n-1}-1$  components. This finishes the proof.
\end{proof}

The following corollary gives an explicit description the $2^{n-1}-1$ components of $\mathcal{V}_{k,n,2}$.

\begin{corollary} \label{cor:rank2}
Each of the $2^{n-1}-1$ irreducible components of
$\,\mathcal{V}_{k,n,2}$ is the affine cone over
a Segre variety $\PP^{ki-1} \times \PP^{k(n-i)-1}$.
There are $\binom{n}{i}$ such components for $i=1,2,\ldots,n-1$.
Hence the reduced scheme of $\,\mathcal{V}_{k,n,2}$ is pure of dimension $kn-1$
and total degree $\,\sum_{i=1}^{n-1} \binom{n}{i} \binom{nk-2}{ik-1}$
\end{corollary}

\begin{proof}
We write our matrix as $\, X = (X_{ij})_{1 \leq i,j \leq n}\,$ where
$X_{ij}$ is a $k \times k$ matrix.
Let $a_1,a_2,\ldots,a_n$ be generic column vectors in $\CC^k$.
For each non-empty subset $I$ of $\{1,2,\ldots,n-1\}$ we set
$$ X_{ij} \,\,= \,\,\begin{cases}\, a_i \cdot a_j^T  &   
\text{ if ($i \in I $ and $j \not\in I$) or  ($i \not\in I $ and $j \in I$)}, \\
                         \,     \quad 0 &
\text{ if $\,\,\, i,j \in I \,$ or  $\,i,j \not\in I $.}
\end{cases}                              $$
Then ${\rm rank}(X)=2$, so this Segre variety $\PP^{k|I|-1} \times \PP^{k(|I^c|)-1}$ sits inside $\mathcal{V}_{k,n,2}$. 
It is  a component since $\mathcal{V}_{k,n,2}$ has dimension $kn-1$. 
The union of these $2^{n-1}-1$ components equals $\mathcal{V}_{k,n,2}$.
\end{proof}

We next turn to the case $r=3$. Here we know that $\mathcal{V}_{k,n,3}$ is
irreducible, but not reduced.

\begin{corollary}
    The variety $\mathcal{V}_{k,n,3}$ has dimension $k(n+1)-3$. It admits the parametrization
    \begin{equation} \label{eq:hadamard}
    X_{ij} \,\,= \,\,(s_i t_j - s_j t_i)^2 \cdot (a_i \cdot a_j^T)  \qquad {\rm for} \,\, 1 \leq i,j \leq n.
    \end{equation}
    This uses $kn+2n$ parameters $a_i \in \CC^k$ and $s_j,t_k \in \CC$. 
    Note that $\mathcal{V}_{k, 3,3}$ is toric for all $k \geq 1$.
\end{corollary}
\begin{proof}
    Using the quadratic form $(x,y,z) \mapsto z^2 - xy$ on $\CC^3$, we can parametrize $n$ vectors on the corresponding quadric in $\CC^3$ by $\begin{bmatrix} s_i^2 & t_i^2 & s_it_i\end{bmatrix}$. 
    The variety $\mathcal{T}_{k,3}$ can then be parametrized by $k \times 3$ matrices of the form $Y_i = a_i \cdot \begin{bmatrix} s_i^2 & t_i^2 & s_it_i\end{bmatrix}$ where $a_i \in \CC^k$ is a column vector. 
    Setting $X_{ij} = Y_iQY_j^T = (s_it_j - s_j t_i)^2 \cdot (a_i \cdot a_j^T)$, we obtain the parametrization~\eqref{eq:hadamard}.
\end{proof}

\begin{example}[$k=2,n=4,r=3$]
The $9$-dimensional variety $\mathcal{V}_{2,4,3}$ is not toric. Its prime ideal
has degree $264$ and it
is generated by $54$ binomial quadrics and $81$ 
six-term quartics.
Each of these is a copy of the
$4 \times 4$ determinant that defines
 $\mathcal{V}_{1,4,3} = {\rm sGr}(2,4)$.
 The number {\bf 4224} in Table \ref{tab:varieties1}
 says that the $4 \times 4$ minors of $X$
define a scheme structure of multiplicity $16$.
\end{example}

We close by examining the cases of interest in physics. Fix $k=2$.
If $r \geq 5$ then the minors generate the prime ideal of $\mathcal{V}_{2,n,r}$.
In our cosmology discussion in Section~\ref{sec3}, 
we had  $r=5$.
An explicit parametrization of 
the irreducible variety $\mathcal{V}_{2,n,r}$ was 
presented in Corollary~\ref{cor:para25}.

The scenario $r=4$ in Section~\ref{sec4} is more interesting.
The variety $\mathcal{V}_{2,n,4}$ has  $2^{n-1}$
irreducible components
$\varphi(\mathcal{Y}_{I,J})$
where $I \sqcup J = [n]$. The main
component $\widetilde{\mathcal{V}}_{2,n,4}$
arises for $I = \emptyset$ and $J = [n]$.
Its parametrization was given in Proposition~\ref{prop:spinorhelicity}.
We now extend this by defining
$$
Y_i^+ \,=\,
\begin{bmatrix}                                                                                               
 u_{i1} v_{i1} & u_{i1} v_{i2} & u_{i2} v_{i2} & \! - u_{i2} v_{i1} \\                                       
\overline{u}_{i1} v_{i1} & \overline{u}_{i1} v_{i2} &                                                        
\overline{u}_{i2} v_{i2} & \! - \overline{u}_{i2} v_{i1} 
\end{bmatrix} 
\quad {\rm and} \quad
Y_i^- \,=\,
\begin{bmatrix}                                                                                              
 u_{i1} v_{i1} & \! - u_{i2} v_{i1} &  u_{i2} v_{i2} & u_{i1} v_{i2} \\
\overline{u}_{i1} v_{i1} & \! - \overline{u}_{i2} v_{i1} &
\overline{u}_{i2} v_{i2} & \overline{u}_{i1} v_{i2} 
\end{bmatrix} .$$

Given any partition $I \sqcup J = [n]$, 
we write $Y_{I,J}$ for the
$2n \times 4$ matrix whose
$i$th block equals $Y_i^+$ if $i \in I$ and equals
$Y_i^-$ if $i \in J$. For $I = \emptyset$, this
 is the matrix we used in
Proposition~~\ref{prop:spinorhelicity}.

\begin{corollary} \label{cor:theirreducible} The irreducible
component $\varphi(\mathcal{Y}_{I,J})$ of $\,\mathcal{V}_{2,n,4}$
has the parametrization
$$  X \,\,\, = \,\,\,
    Y_{I,J} \cdot \begin{small}
    \begin{bmatrix}
        0 & \! \rm{Id}_2 \\ \rm{Id}_2 &\! 0
    \end{bmatrix} \end{small} \cdot (Y_{I,J})^T .$$
\end{corollary}

\begin{proof}
The matrices $Y_i^+$ and $Y_i^-$ give the
components
$\mathcal{T}^+_{2,4}$
and $\mathcal{T}^-_{2,4}$
of the variety~$\mathcal{T}_{2,4}$.
\end{proof}

Corollary \ref{cor:theirreducible} is used to study degrees and equations of the $2^{n-1}$ components of $\mathcal{V}_{2,n,4}$.

\begin{example}[$n = 4$]
The $14$-dimensional variety $\mathcal{V}_{2,4,4}$ has
degree {\bf 2772}; see Table~\ref{tab:varieties1}.
Its eight components are indexed by $I = 
\emptyset,\{1\},\{2\},\{3\},\{4\},
\{1,2\},\{1,3\},\{1,4\}$.
Their degrees are $412,344,344,344,344,328,328,328$,
for a total of $2772$.
The prime ideals of the irreducible components of $\mathcal{V}_{2,4,4}$ are generated by polynomials of
degree $2,3,4,5$.
The quadrics are $2 \times 2$ minors and
the quintics are $5 \times 5$ minors.
Table \ref{tab:V244}
gives the number of minimal generators.

\begin{center}
\begin{table}[h]
\begin{center}
    \begin{tabular}{ccccccc}
     $I$                              && \text{degree} & {\# quadrics} & {\# cubics} & {\# quartics} & {\# quintics} \\
       $\emptyset$                    && 412  & 0  &  60 & 12 & 108 \\ 
       $\{1\}, \{2\}, \{3\}, \{4\}$   && 344  & 3  & 45   & 30  & 108 \\ 
       $\{1,2\}, \{1,3\}, \{1,4\}$    && 328  & 4  & 40  & 36  & 108
\end{tabular}
\end{center}
\caption{Degree and minimal generators of the prime ideals of the 8 components in $\mathcal{V}_{2,4,4}$. \label{tab:V244} }
\end{table}    
\end{center}
\end{example}

\bigskip

\noindent {\bf Acknowledgment}:
We thank Ben Hollering and Chia-Kai Kuo for help with this project.

\bigskip

\begin{scriptsize}
\noindent {\bf Funding statement}:
Research supported by the European Union (ERC, UNIVERSE PLUS, 101118787). Views and opinions 

\noindent expressed are however those of the authors only and do not necessarily reflect those of the European Union or the European

\noindent Research Council Executive Agency. Neither the European Union nor the granting authority can be held responsible for them.

\vspace{2mm}

\noindent {\bf Data availability statement}: The code used to compute the results in this paper is available upon request from the corresponding author.

\end{scriptsize}

\bigskip
\bigskip

\footnotesize{
{\bf Authors' addresses:}

\bigskip

Yassine El Maazouz, Caltech \hfill {\tt \href{mailto:maazouz@caltech.edu}{maazouz@caltech.edu}}

Bernd Sturmfels, MPI MiS Leipzig  \hfill {\tt \href{mailto:bernd@mis.mpg.de}{bernd@mis.mpg.de}}}

Svala Sverrisd\'ottir, UC Berkeley  \hfill {\tt \href{mailto:svala@math.berkeley.edu}{svala@math.berkeley.edu}}

\end{document}